\documentclass{amsart}
\usepackage{amssymb,latexsym,graphicx}

\DeclareMathOperator{\supp}{supp}

\DeclareMathOperator{\var}{var}

\def\L{\mathcal{L}}
\def\J{\mathcal{J}}

\def\bv{{\mathcal BV}}

\def\R{\mathbb R}
\def\Z{\mathbb Z}
\def\N{\mathbb N}
\def\C{\mathbb C}

\def\zeroOne{zero-one }
\def\Q{{\bf Q}}
\def\pwm{piecewise monotonic}
\def\htau{{\widehat\tau}}

\def\mapright#1{\smash{\mathop{\longrightarrow}\limits^{#1}}}

\def\mapright#1{\smash{\mathop{\longrightarrow}\limits^{#1}}}

\theoremstyle{plain}

\newtheorem{theorem}{Theorem}[section]
\newtheorem{proposition}[theorem]{Proposition}
\newtheorem{lemma}[theorem]{Lemma}
\newtheorem{corollary}[theorem]{Corollary}

\theoremstyle{definition}
\newtheorem{definition}[theorem]{Definition}
\newtheorem{example}[theorem]{Example}

\theoremstyle{remark}
\newtheorem*{notation}{Notation}

\newif\ifproofing
\proofingfalse
\newcommand{\note}[1]{\ifproofing\marginpar{#1}\fi}

\newcommand{\prop}[2]{\begin{proposition}\label{#1}\note{#1}#2\end{proposition}}
\newcommand{\lem}[2]{\begin{lemma}\label{#1}\note{#1}#2\end{lemma}}
\newcommand{\cor}[2]{\begin{corollary}\label{#1}\note{#1}#2\end{corollary}}
\newcommand{\theo}[2]{\begin{theorem}\label{#1}\note{#1}#2\end{theorem}}
\newcommand{\defi}[2]{\begin{definition}\label{#1}\note{#1}#2\end{definition}}
\newcommand{\exe}[2]{\begin{example}\label{#1}\note{#1}#2\end{example}}
\newcommand{\notate}[1]{\begin{notation}#1\end{notation}}

\newcommand{\prooff}[1]{\begin{proof}#1\end{proof}}

\begin{document}

\renewcommand{\theenumi}{\roman{enumi}}
\renewcommand{\labelenumi}{(\theenumi)} 

\title{Dimension groups for interval maps II:\\ The transitive case}
\author{Fred Shultz}
\address{Wellesley College\\Wellesley, Massachusetts 02481}
\email{fshultz@wellesley.edu}

\begin{abstract}

Any continuous,  transitive, \pwm\ map is determined up to a binary choice by its  
dimension module with the associated finite sequence of generators. 
The dimension module by itself  determines the topological entropy
of any transitive \pwm\ map, and determines any transitive unimodal map up to conjugacy. 
For a transitive \pwm\ map which is not essentially injective, the associated dimension group
is a direct sum of simple dimension groups, each with a unique state.

\end{abstract}

\maketitle

\section{Introduction}

Motivated by Elliott's \cite{Ell} classification of AF-algebras by dimension groups, Krieger
\cite{KriShift} gave a dynamical definition of a dimension group for ample groups of
homeomorphisms on zero dimensional compact metric spaces in terms of an equivalence relation on
the compact open subsets, and used this to associate  dimension groups with (two-sided) shifts of
finite type. He showed that one of these dimension groups (together with an associated automorphism) 
is a complete invariant for
shift equivalence of topologically mixing  shifts of finite type.

The current paper is part of a program to define and investigate dimension groups for piecewise
monotonic maps of the unit interval, with the goal of obtaining interesting invariants for such
maps. We will see that there are good analogues for interval maps of
several well known results for shifts of finite type.

In \cite{Shu}, a dimension group is defined for \pwm\ maps of the unit interval. The basic
procedure is to associate with the given \pwm\ map $\tau:I\to I$ a closely related local
homeomorphism $\sigma:X\to X$, where $X$ is formed by disconnecting the interval at
appropriate points.  Then, following Boyle, Fiebig, and Fiebig \cite{BFF}, the dimension group
$DG(\tau)$ is defined to be $C(X,\Z)$, modulo the equivalence relation given by  $f \sim g$ if
$\L^n f = \L^n g$ for some $n \ge 0$, where $\L= \L_\sigma$ is the transfer operator
$$(\L_\sigma f)(x) = \sum_{\sigma y = x} f(y).$$
This
dimension group comes equipped with a natural injective endomorphism $\L_*:[f] \mapsto [\L
f]$, which is an automorphism if $\tau$ is surjective. The triple ($DG(\tau)$, $DG(\tau)^+$, 
$\L_*$) is called the dimension triple associated with $\tau$.
In the current paper, we strengthen the results in \cite{Shu} by restricting consideration to
transitive maps. 

For a transitive one-sided shift of finite type
$(X_A,\sigma_A)$, there is a period
$p$, and a partition of $X$
 into $p$ clopen pieces permuted by $\sigma$, such that
$\sigma^p$ is topologically mixing on each. The associated dimension group for
$\sigma_A$ will be a direct sum of simple dimension groups.

We establish similar results for transitive \pwm\ maps. If $\tau$ is
transitive, and the associated local homeomorphism $\sigma$ is not a homeomorphism, then there
exists a positive integer
$N$ and a decomposition of
$X$ into clopen sets
$X_1, \ldots, X_N$, permuted by $\sigma$, such that $\sigma^N$ is topologically exact on each
$X_i$ (Theorem \ref{1.66}).  (A map $\sigma:X_i\to X_i$ is topologically exact if for every
non-empty open set $V\subset X_i$, there exists $n \ge 0$ such that $\sigma^n(V) = X_i$.) The dimension
group $DG(\tau)$ will be a direct sum of simple dimension groups (Theorem \ref{1.73}).
(Ironically, the case where $\sigma$ is bijective is the bad case: such a
decomposition need not exist for such maps.)

The same kind of decomposition
holds for $\tau$, except that the clopen sets are replaced by finite unions of closed intervals,
which may overlap at endpoints. This result for $\tau$ (with topologically mixing instead of topologically exact)  can be
obtained from the decomposition of the non-wandering sets of
\pwm\ maps due to Hofbauer \cite[Thm. 2]{HofDecomp}. General results on ``regular periodic decompositions" of
(continuous) transitive maps can be found in \cite{Banks}.

The order on a dimension group is an important part of its structure, and for transitive maps
we can give an effective description of that order. For a simple dimension group, 
the order is determined by the states, i.e., the normalized order
preserving homomorphisms into $\R$, cf., e.g., \cite{Eff}. If $\tau:I\to I$ is topologically
exact, the dimension group $DG(\tau)$ is simple, with a unique state (Corollary
\ref{1.75}).

More generally, when $\tau$ is transitive, there is a  state that is scaled by $\L_*$ by a
factor $s \ge 1$. This state is given by a measure on $I$ that is scaled by $\tau$ by
the factor $s$. If $\tau$ is essentially injective (i.e., there are no intervals $J_1$, $J_2$
with $\tau(J_1) = \tau(J_2)$), or equivalently, if the associated local homeomorphism $\sigma$
is not a homeomorphism, then transitivity of $\tau$ implies that the state scaled by $\L_*$ is
unique, given by the unique measure on $I$ scaled by $\tau$, with the scaling factor $s = \exp h_\tau$, where $h_\tau$ is the
topological entropy (Theorem \ref{1.73}).  It follows that for transitive maps, the dimension triple determines the
topological entropy of the map. Transitive unimodal maps are determined by their topological
entropy, and thus two such maps are determined up to conjugacy by their dimension triples.

If $\tau$ is surjective, then $\L_*$ is an automorphism of $DG(\tau)$, so the latter can be
viewed as a $\Z[t,t^{-1}]$ module, called the dimension module. In
\cite{Shu}, a canonical finite sequence of generators for this module is identified.  The dimension
module with its canonical  sequence of generators, determines a continuous transitive
 map up to a binary choice, namely, whether the function increases or decreases on the first
interval of monotonicity (Theorem \ref{1.79}).

For a two-sided irreducible shift of finite type, Krieger \cite{KriShift} showed that the
measure of maximal entropy gives an imbedding of the dimension group modulo the subgroup of
infinitesimals into $\R^n$, and that the range of this imbedding on the clopen subsets is an  easily
computable invariant,  powerful enough to distinguish mixing shifts with the same zeta function. This quotient group has
also proved to be a valuable invariant for minimal homeomorphisms of Cantor sets, cf. \cite{GPS}. 

If $\tau:I\to I$ is a transitive \pwm\ map, the infinitesimal
elements  in
$DG(\tau)$ are precisely those killed by all states. We show that the unique measure on the
unit interval
$I$ scaled by $\tau$ induces an isomorphism from $DG(\tau)/DG(\tau)_{\rm inf}$ into $\R^n$
(Proposition \ref{1.77}). This provides an alternative invariant, which is  easy to
compute in terms of the measure scaled by $\tau$.

For a primitive matrix $A$, the Perron-Frobenius theorem establishes
existence of a unique eigenvector for the maximum eigenvalue
$\lambda$, and convergence of powers of $\lambda^{-1}A$.  In the current paper, the role of
$\lambda^{-1}A$ is played by $s^{-1}\L_\sigma*$, which is  the Perron-Frobenius
operator.  Convergence of powers of this operator on the space $\bv$ of functions in $L^1$ of
bounded variation has been thoroughly investigated, e.g. \cite{RuStat, HofKel, Ry, WalBeta, WalExp, RuBook}. An
exposition can be found in the book \cite{BG}, and we mainly follow
\cite{Ry}.
For our purposes, it is important to have not just convergence in the $\bv$ norm, but uniform
convergence of continuous functions.  Since $\sigma$ is a local homeomorphism, then $C(X)$ is
invariant under $\L_\sigma$, while $C(I)$ is not invariant under $\L_\tau$. Thus we work with $\sigma$ rather than
with
$\tau$, and establish the necessary results on continuity and uniform convergence in an appendix.  Similar uniform
convergence results for the case of
$\beta$-transformations were established by Walters \cite{WalBeta}. These results are the key technical tools used to
establish the decomposition of transitive maps into topologically exact pieces, and the results on uniqueness of
scaled states described above.

We now summarize this paper.  We begin with background on dimension groups of interval maps, 
mostly taken from \cite{Shu}. Then scaling measures are introduced, 
which are measures $\mu$ on $I$ such that there is $s > 0$ such that
$\mu(\tau(E)) = s \mu(E)$ for all Borel sets $E$ on which $\tau$ is injective.  These are a special case of
conformal measures (\cite{Ur-De}), and are closely related to conjugacies of \pwm\ maps with uniformly piecewise
linear maps, i.e., piecewise linear maps with slopes $\pm s$.  (By a result of Parry \cite{Parry}, every transitive
\pwm\ map is conjugate to a uniformly piecewise linear map, cf. Corollary \ref{1.65}.) This is mostly standard
material, adapted to the current context.  For our purpose, the key fact is that such measures induce states on the
dimension group.  

We next prove the finite decomposition of non-injective transitive maps into topologically
exact maps, and the uniqueness of the state scaled by $\L_*$. A similar uniqueness result for topologically exact
positively expansive maps was proven by Renault \cite{RenRN}, who used a somewhat different dimension group, based
on inductive limits of $C(X)$ rather than $C(X,\Z)$.

 We finish by computing the dimension triple (or the
quotient by the subgroup of infinitesimals) for several families of uniformly  piecewise linear maps, including tent maps and
$\beta$-transformations. In the latter case, Katayama, and Watatani
\cite{KatMatWat} have associated C*-algebras $F_\beta^\infty$ and $O_\beta$ with the
$\beta$-transformation, and $K_0(F_\beta^\infty)$ is a dimension group. We show that $K_0(F_\beta^\infty)$ is isomorphic
to  $DG(\tau)$ as a group (and is order isomorphic if the orbit of 1 is eventually periodic.) 

The author would like to thank David Handelman for fruitful conversations regarding this paper, and in particular for pointing
out the role of convergence of powers of the Perron-Frobenius operator in proving uniqueness of states.

\section{Background}

We first review the construction in \cite{Shu}, which to each piecewise monotonic map
$\tau:[0,1]\to [0,1]$ associates a local homeomorphism $\sigma:X\to X$ on a compact set $X\subset \R$, where $X$
is constructed by disconnecting
$[0,1]$ at certain points.  (For references to related work of other authors, see \cite{Shu}.)

\defi{1.1}{Let
$I = [0,1]$. A map
$\tau:I
\to I$ is  \emph{piecewise monotonic} if there are points
$0 = a_0 < a_1 <
\ldots < a_n = 1$ such that $\tau|(a_{i-1},a_i)$ is continuous and strictly monotonic for $1 \le i \le n$. 
We will assume the sequence $a_0, a_1, \ldots, a_n$ is chosen so that no interval
$(a_{i-1},a_i)$ is contained in a larger open interval on which
$\tau$ is continuous and strictly monotonic. The sequence of points $0 = a_0 < a_1 < \cdots < a_n =
1$ is the {\it partition associated with $\tau$}, and the intervals
$\{(a_{i-1}, a_i)\mid 1\le i\le n\}$ are called the \emph{intervals of monotonicity for $\tau$}. Note
that for
$1
\le i
\le n$, the map
$\tau|(a_{i-1},a_i)$ extends uniquely to a strictly monotonic continuous map $\tau_i:[a_{i-1},a_i] \to I$, which will be
a homeomorphism onto its image.}

If $\tau$ is not continuous,  we will ignore the actual values of $\tau$ at the partition
points, and instead view $\tau$ as being multivalued, with the values at $a_i$ (for $1 \le i \le n-1$) being the
values given by left and right limits, i.e. the values of $\tau$ at $a_i$ are $\tau_i(a_i)$ and
$\tau_{i+1}(a_i)$. We define a (possibly multivalued) function $\htau$ on $I$ by setting $\htau(x)$ to
be the set of left and right limits of $\tau$ at $x$.  At points where $\tau$ is continuous, $\htau(x) = \{\tau(x)\}$, and we
identify $\htau(x)$ with $\tau(x)$.
Thus for $A\subset I$, $\htau(A) = \bigcup_i \tau_i(A\cap [a_{i-1},a_i])$, and $\htau^{-1}(A) =
\bigcup_i \{x \in [a_{i-1},a_i] \mid \tau_i(x) \in A\}$.

If $x \in I$,  the \emph{generalized
orbit} of $x$ is the smallest subset of $I$ containing $x$ and forward and backward invariant with respect to
$\htau$. Let
$I_1$ be the union of the generalized orbits of $a_0, a_1,  \ldots, a_n$, and let
$I_0 = I \setminus I_1$.

\defi{1.1.1}{Let $I = [0,1]$, and  let $I_0$, $I_1$ be as above. 
The {\it disconnection of $I$ at points in
$I_1$\/} is the totally ordered set $X$ which consists of  a copy of $I$ with the usual ordering,  but with each point
$x
\in I_1\setminus
\{0,1\}$ replaced by two points $x^- < x^+$. We equip $X$ with the order topology, and define the {\it collapse map\/}
$\pi:X\to I$ by
$\pi(x^\pm) = x$ for $x \in I_1$, and $\pi(x) = x$ for $x \in I_0$. We write  $X_1 = \pi^{-1}(I_1)$, and $X_0 = \pi^{-1}(I_0)
=X\setminus X_1$.}

\notate{For any pair $a, b \in
X$, we write
$[a,b]_X$ for the order interval
$\{x
\in X
\mid a
\le x
\le  b\}$.  If $b_1, b_2 \in I_1$ with $b_1 <
b_2$, then
$I(b_1, b_2)$ is the  order interval
$[b_1^+,b_2^-]_X$. If $b_1 > b_2$, then we define $I(b_1,b_2) = I(b_2,b_1)$, and if $b_1 = b_2$, then we set $I(b_1,b_2) =
\emptyset$.}

\prop{1.2}{(\cite[Prop. 2.2]{Shu}) Let $X$ be the disconnection of $I$ at points in $I_1$, and $\pi:X\to I$ the collapse
map. Then
$X$ is homeomorphic to a compact subset of
$\R$, and
\begin{enumerate}
\item $\pi$ is  continuous and order preserving.
\item $I_0$ is dense in $I$, and $X_0$ is dense in
$X$.
\item $\pi|_{X_0}$ is a homeo\-morphism from $X_0$ onto $I_0$.
\item $X$ has no isolated points.
\item If $a, b \in I_1$, then $I(a,b)$ is clopen in $X$, and every clopen subset of $X$ is a
finite disjoint union of such order intervals.
\end{enumerate}}

\prop{1.2.1}{(\cite[Thm. 2.3]{Shu}) Let $\tau:I\to I$ be a piecewise
monotonic map, with associated partition $a_0 < a_1 < \cdots < a_n$, and let $(X,\pi)$ be as described in
Proposition \ref{1.2}.
\begin{enumerate}
\item  There is a unique continuous map
$\sigma:X\to X$ such that $\pi\circ \sigma = \tau \circ \pi$ on $X_0$.\label{I:1.2.1.1}
\item The sets $X_0$ and $X_1 = X\setminus X_0$ are forward and backward invariant with respect to $\sigma$.\label{I:1.2.1.1.5}
\item $\pi$ is a conjugacy from $\sigma|X_0$ onto $\tau|I_0$.\label{I:1.2.1.2}
\item The sets $J_1 = I(a_0,a_1), \ldots, J_n= I(a_{n-1},a_n)$  are a partition of $X$ into clopen sets such
that for $1\le i\le n$, $\pi(J_i) = [a_{i-1},a_i]$, and  $\sigma|J_i$ is a homeomorphism from $J_i$ onto
the clopen set $\sigma(J_i)$.\label{I:1.2.1.3}
\item For $1\le i \le n$, $\pi\circ \sigma = \tau_i\circ \pi$ on $J_i$.\label{I:1.2.1.4}
\end{enumerate}
}

Note that if $\tau$ is continuous, by (\ref{I:1.2.1.2})  and density of $I_0$ in $I$ and $X_0$ in $X$,  $\pi$ will be a
semi-conjugacy from $(X,\sigma)$ onto $(I,\tau)$.

If $\tau:I\to I$ is \pwm, and $\sigma:X\to X$ is the map in Proposition \ref{1.2.1}, then $\sigma$
will be a local homeomorphism by Proposition \ref{1.2.1} (\ref{I:1.2.1.3}). 
We will call $\sigma$ the \emph{local homeomorphism 
associated with $\tau$}.   Property (\ref{I:1.2.1.2})  of Proposition
\ref{1.2.1} can be used to show that $\sigma$ and $\tau$ share many properties. Before being more explicit, we review some
terminology.

\defi{1.7}{If $X$ is any topological space, and $f:X\to X$ is
a continuous map, then
$f$ is \emph{transitive} if for each pair
$U, V$ of non-empty open sets, there exists $n \ge 0$ such that
$f^n(U)
\cap V
\not=
\emptyset$. We say $f$ is \emph{strongly transitive} if for
every non-empty open set $U$, there exists $n$ such that
$\cup_{k=0}^nf^k(U) = X$, and $f$ is \emph{topologically exact} if for
every non-empty open set $U$, there exists $n$ such that $f^n(U)= X$.}

The notion of strong transitivity appears in Parry's paper \cite{Parry}, and a discussion of topological exactness
for interval maps can be found in the book of Preston, cf., e.g.,  \cite[pp. 6, 27]{Preston}.

\defi{1.7.1}{
 If $\tau:I\to I$ is \pwm, we view $\tau$ as undefined at the set $C$ 
of endpoints of intervals of monotonicity, and say $\tau$ is  \emph{transitive}
if  for each pair
$U, V$ of non-empty open sets, there exists $n \ge 0$ such that
$\tau^n(U)
\cap V
\not=
\emptyset$. We say $\tau$ is {\it strongly transitive\/} if for
every non-empty open set $U$, there exists $n$ such that $\cup_{k=0}^n\htau^k(U) = I$.
(Recall that
$\htau$ denotes the (possibly multivalued) function whose value
at each point $x$ is given by 
 the left and right hand limits of $\tau$ at $x$.) The map $\tau$ is \emph{topologically exact} if 
for
every non-empty open set $U$, there exists $n$ such that $\htau^n(U) = I$. (See Lemma \ref{1.89}, Example \ref{1.98}, and Lemma
\ref{1.105} for examples of topologically exact maps.)}

If $\tau:I\to I$ is continuous and \pwm, both Definitions \ref{1.7} and \ref{1.7.1} are applicable, and are consistent.
Transitivity is equivalent to the existence of a dense orbit.
 
The full two sided $n$-shift is an example of a map that is transitive, but not strongly transitive. (The
complement of a fixed point is an invariant open set.) However, for
\pwm\ maps, these notions are equivalent.

\prop{1.3}{Let $\tau:I\to I$ be \pwm, and  $\sigma:X\to X$ the associated local homeomorphism.
\begin{enumerate}
\item $\sigma$ is surjective iff $\tau$ is surjective.\label{I:1.3.1} 
\item $\sigma$ is strongly transitive iff $\tau$ is strongly transitive iff $\tau$ is
transitive iff $\sigma$ is  transitive.\label{I:1.3.2} 
\item $\sigma$ is topologically exact iff $\tau$ is topologically exact.\label{I:1.3.3} 
\end{enumerate}}

\prooff{\cite[ Lemma 4.2, Prop. 2.9, Lemma 5.2]{Shu}.}

Now we turn to comparing the topological entropies of $\tau$ and $\sigma$. For a
continuous map $f$ on a compact metric space
$X$, we denote the (topological) entropy of $f$ by $h_f$. In the special case where $(X,f)$ is a
subshift, then we have the  formula
\begin{equation}
h_f = \lim_{n\to \infty} \textstyle{\frac{1}{n}} \ln c_n(f),
\end{equation}
where $c_n(f)$ is the number of cylinders in $X$ of length $n$, cf. \cite[ex. 6.3.4]{Lind-M}.

If $\tau:I\to I$ is \pwm, let  $a_0 < a_1 < \ldots < a_q$ be the  partition associated with
$\tau$, and let $A_j = (a_{j-1},a_j)$ for $1\le j\le q$. Let 
\begin{equation}
\Sigma_\tau = \left\{(s_i)_{i=0}^\infty \in \prod_0^\infty\{1,
\ldots, q\} \mid \bigcap_{i=0}^n  \tau^{-i}(A_{s_i})\not= \emptyset \hbox{ for all $n \ge 0$}\right\}.
\end{equation}
Then $\Sigma_\tau$ is compact (for the product topology), and invariant under the shift map. The
shift map on $\Sigma_\tau$  is called  the symbolic dynamics of $\tau$. If $\tau$ is continuous as
well as \pwm, then the topological entropy of $\tau$ is the same as the topological entropy of the
symbolic dynamics of $\tau$, cf. \cite[Prop. 2.1]{Young} or \cite{Mis-Zie}.  If $\tau$ is \pwm, but is not
continuous everywhere, then we define $h_\tau$ to be the topological entropy of the symbolic dynamics
of $\tau$.

\prop{1.9}{If $\tau:I\to I$ is \pwm, and $\sigma:X\to X$ is the associated local
homeomorphism, then $h_\tau = h_\sigma$.}  

\prooff{It is readily seen that the symbolic dynamics of $\tau$ and of $\sigma$ coincide (where the
symbolic dynamics of $\sigma$ is taken with respect to the partition $J_1, \ldots,
 J_n$ in Theorem \ref{1.2.1}.) As remarked above, by \cite[Prop. 2.1]{Young} the topological entropy of
$\tau$ is the same as that of its symbolic dynamics.  The proof of \cite[Prop. 2.1]{Young} is given
for a continuous \pwm\ map
$\tau:I\to I$, but applies without change to
$\sigma:X\to X$. Thus  $\sigma$ and its symbolic dynamics have the same topological
entropy, and so $h_\tau = h_\sigma$.}

\defi{1.22}{Let $X$ be a compact Hausdorff space and $\sigma:X\to X$ any map such that all fibers $\sigma^{-1}(x)$ are 
finite.  Then for any $f:X\to \R$, we define $\L_\sigma f$ by 
$$(\L_\sigma f)(x) = \sum_{\sigma y = x} f(y),$$
and call $\L_\sigma$ the \emph{transfer map}.
We will write $\L$ in place of $\L_\sigma$ when the meaning is clear from the context.}

If $\sigma:X\to X$ is a local homeomorphism, and $f:X\to \R$ is continuous, then $\L_\sigma f$ will be continuous (see the
remark after Lemma \ref{1.55}), so
$\L_\sigma$ maps $C(X)$ into $C(X)$, and $C(X,\Z)$ into $C(X,\Z)$.

 If
$X$ is a compact metric space, a map
$\sigma:X\to X$ is a
\emph{piecewise homeo\-morphism} if
$\sigma$ is continuous and open, and $X$ admits a finite partition into clopen sets $X_1, X_2,
\ldots, X_n$ such that $\sigma$ is a homeo\-morphism from $X_i$ onto $\sigma(X_i)$ for $i = 1, \ldots, n$. Any local homeomorphism
on a zero dimensional compact metric space will be a piecewise homeomorphism, and if $\tau:I\to I$ is \pwm, then the associated
local homeomorphism $\sigma:X\to X$ will be a piecewise homeomorphism (Proposition \ref{1.2.1}).

\defi{1.3.1}{Let $X$ be a compact metric space, and $\sigma:X\to X$ a piecewise homeomorphism. For $f, g \in C(X,\Z)$, 
define $f\sim g$ if there exists $n \ge 0$ such that $\L^n f = \L^n g$. Then $G_\sigma$ is the
ordered abelian group whose elements are equivalence classes $[f]$, with addition $[f] + [g] = [f+g]$, and with order
given by
$[f]\ge 0$ if
$\L^n f\ge 0$ for some
$n$. If $\tau:I\to I$ is \pwm, with associated local homeomorphism $\sigma:X\to X$, then $DG(\tau)$ is defined to be
$G_\sigma$.}

The ordered groups $G_\sigma$ and $DG(\tau)$ are dimension groups \cite[Cor. 3.12 and Def. 3.13]{Shu}, i.e., are  inductive
limits of a sequence of groups of the form $\Z^{n_k}$.  See \cite{Eff} or \cite{Goo} for background on dimension groups.

\defi{1.16}{$\L_*:DG(\tau)
\to DG(\tau)$ is defined by $\L_*[f] = [\L f]$.}

This is an injective  homomorphism, and is bipositive, i.e., $[f] \ge 0$ iff $\L_*[f] \ge 0$. If $\tau$ is surjective, then
$\L_*$ is surjective, and thus is an automorphism of the dimension group $DG(\tau)$.  In that case, we define an action of
$\Z[t,t^{-1}]$ on
$DG(\tau)$ by
$$(\sum_{-n}^n a_it^i)[f] = (\sum_{-n}^n a_i\L_*^i)[f],$$ and view $DG(\tau)$ as a $\Z[t,t^{-1}]$ module. The following
describes a  set of generators for this module. In the statement of the  theorem, we identify order
intervals $I(a,b)$ in $X$ with the equivalence class in
$DG(\tau)$ of their characteristic functions.

\theo{1.41}{Let $\tau:I\to I$ be \pwm\ and  surjective, with associated partition $C = \{a_0, a_1,
\ldots, a_n\}$.  
Let
\begin{equation}
\J_1 =
\{\,I(c,d)
\mid \hbox{$c$, $d$ are adjacent points in 
$\{\,a_0, a_1,
\ldots, a_n\,\}$}\,\},
\end{equation}
and let $\J_2$ be the set of intervals corresponding to jumps at partition points, i.e., 
\begin{equation}
\J_2 = \{\,I(\tau_i(a_i), \tau_{i+1}(a_i)) \mid 1 \le i
\le n-1\,\}.
\end{equation}
 Then
$DG(\tau)$ is generated as a module by $\J_1 \cup \J_2$.}

\prooff{\cite[Thm. 6.2]{Shu}}

\cor{1.42}{If  $\tau:I\to I$ is a
continuous, surjective piecewise monotonic map with associated partition
$\{a_0, a_1, \ldots, a_n\}$, then
$DG(\tau)$ is generated as a module by $I(a_0,a_1)$, $I(a_1,a_2)$, $\ldots$,
$I(a_{n-1},a_n)$.}

\prooff{\cite[Cor. 6.3]{Shu}}

\section{Scaling measures, and uniformly piecewise linear maps}

For simple dimension groups, the order is
determined by states, i.e., positive homomorphisms into $\R$, cf. \cite{Eff}. 
In the case that $\tau$ is transitive, we will see that  states are given by scaling measures. 
In this
section, we provide background for such measures, and discuss
their connection with uniformly piecewise linear maps.

Let $X$ be a compact metric space and $m$ a probability measure on $X$ (i.e. a
positive regular Borel measure with $m(X) = 1$.) If
$f
\in L^1(X,m)$, we will usually write
$m(f)$ instead of $\int f \,dm$.

\defi{1.54}{Let $X$ be a compact metric space and $m$ a probability measure on $X$, and 
$\sigma:X\to X$ a  map that takes Borel sets to Borel sets.
Then {\it $\sigma$ scales $m$ by a factor $s$} if
$m(\sigma(E)) = s\,m(E)$ for all Borel sets $E$ on which
$\sigma$ is 1--1. (This is a special case of the notion of a {\it conformal measure\/}, cf. \cite{Ur-De}.)
 The measure $m$ has {\it
full support\/} if its support is all of
$X$.}

Now we will show that scaling measures are just the eigenvectors of the dual of the
transfer operator.  We write $C(X)$ for the Banach space of real valued continuous
functions with the supremum norm.  Recall for a local homeomorphism $\sigma:X\to
X$, 
$\L$ is the transfer operator
$\L_\sigma$ (Definition \ref{1.22}).

The following result, in the more general context of conformal measures, can be found in \cite[Prop. 2.2]{Ur-De}.

\lem{1.55}{Let $\tau:I\to I$ be \pwm, and let $\sigma:X\to X$ be the associated local homeomorphism.  Let
$\mu$ be a probability measure on $X$, and
$0 < s \in \R$. The following are equivalent.
\begin{enumerate}
\item $\mu(\L f) =
s\mu(f)$ for all $f \in C(X)$.
\item   $\sigma$ scales $\mu$ by the factor $s$.
\item $\L$ is well defined on $L^1(X,\mu)$, and $\mu(\L f) =
s\mu(f)$ for all $f \in L^1(X,\mu)$.
\end{enumerate}}

\prooff{(i) $\Rightarrow$ (ii). Let $J_1, \ldots, J_n$ be a partition of $X$ into clopen sets such that $\sigma$
is injective on each $J_i$, cf. Proposition \ref{1.2.1}. Let $E_i = \sigma(J_i)$, and let $\psi_i:E_i\to J_i$ be the
inverse of
$\sigma|J_i$. Let $f\in C(X)$ have support in $J_i$. Then $\L f$ has support in
$E_i$ and is defined on $E_i$ by $\L f = f\circ \psi_i$. Thus
\begin{equation}
s\int_{J_i}\, f\,d\mu = s\mu(f) = \mu(\L f) = \int_{E_i} f\circ \psi_i \, d\mu =
\int_{J_i} f \, d(\psi_i\cdot \mu),\label{(1.29)}
\end{equation}
where $(\psi_i\cdot \mu)(E) = \mu(\sigma(E))$ for $E \subset J_i$. Since this holds
for all $f \in C(J_i)$, it follows that the regular Borel measures $s\mu$ and
$\psi_i\cdot \mu$ must coincide on $J_i$. Thus for $E \subset J_i$ we have
$\mu(\sigma(E)) = s\mu(E)$. Therefore $\mu(\sigma(E\cap J_i)) = s\mu(E\cap J_i)$
for all Borel sets $E$. It follows that $\mu(\sigma(E)) = s\mu(E)$ for all Borel
sets $E$ on which $\sigma$ is injective. 

(ii) $\Rightarrow$ (iii). If $A$ is a Borel set and 
$\mu(A) = 0$, then $\mu(\sigma(A\cap J_i)) = 0$ for all $i$, so $\mu(\sigma(A)) = 0$.
If $f = 0 $ a.e., since 
$$\{x \mid (\L f)(x) \not= 0\} \subset \sigma(\{x \mid f(x) \not= 0\}),$$
then $\L f = 0$ a.e.
 Thus
$\L$ is well-defined on
$L^1(X,\mu)$. Reversing the argument in the first paragraph (but with $f\in L^1(X,\mu)$ in (\ref{(1.29)})) proves (iii).
Trivially (iii) implies (i).}

We note for future reference that the proof of Lemma \ref{1.55} shows that $\L_\sigma$ maps $C(X)$ into $C(X)$.
\smallskip

Recall that a probability measure $m$ on a compact metric space is {\it non-atomic\/} if $m$
takes the value $0$ on each singleton subset $\{x\}$ of $X$.

\prop{1.56}{Let $\tau:I\to I$ be \pwm, with associated local
homeomorphism
$\sigma:X\to X$, and collapse map $\pi:X\to I$. Then $\mu \mapsto \mu\circ \pi$ is a 1-1 correspondence of
non-atomic probability measures  on $I$ of full support, scaled by $\tau$ by a factor $s$, and
non-atomic probability  measures of full support on
$X$, scaled by
$\sigma$ by the factor $s$.}

\prooff{Since $I_1$ and $X_1$ are countable, these sets are killed by any non-atomic
measures. The proposition then follows from the fact that the collapse map $\pi$ is a
conjugacy from $\sigma|X_0$ onto $\tau|I_0$ (Proposition \ref{1.2.1}).}

\defi{1.57}{A map $f:I\to I$ is {\it piecewise linear\/} if $f$ is \pwm\ and is linear on the
interior of each interval of monotonicity. A map is {\it uniformly piecewise linear\/} if it is piecewise linear and the
slopes are all $\pm s$ for some $s> 0$.}

\exe{1.58}{If $f:I\to I$ is uniformly piecewise linear with slopes $\pm s$, then Lebesgue
measure is scaled by $f$ by a factor $s$.}

This example is canonical, in the sense made precise by   Proposition \ref{1.59}. This result is implicit
in \cite{Parry}.

\prop{1.59}{A \pwm\ map $\tau:I\to I$  is
conjugate to a piecewise linear map with slopes $\pm s$ iff there exists a non-atomic
probability measure $m$ on $I$ of full support on $I$, scaled by $\tau$ by the factor
$s$. The conjugacy can be chosen to carry $m$ to Lebesgue measure.}

\prooff{Suppose such a measure $m$ exists. Define $h:I \to I$ by $h(x) = m((0,x])$ (with
$h(0) = 0$).  Since $m$ is zero on no open set,  then $h$ is
strictly increasing. The map $h$ is right continuous by $\sigma$-additivity of
$m$, and left continuous  because
$m$ is non-atomic. Thus $h$ is a continuous injective map of $[0,1]$ into
$[0,1]$. Since its range includes 0 and 1 and is connected, then $h$ is 
surjective, so is a homeomorphism. 

Now define $f:[0,1] \to
[0,1]$ by $f= h \circ\tau\circ h^{-1}$. Then $f:[0,1]\to [0,1]$  is 
conjugate to $\tau$.   Let $J_1,
\ldots, J_n$ be a cover of
$I$ by intervals on which
$\tau$ is monotonic, with only endpoints in common.  From
the fact that $m$ is scaled by $\tau$ by a factor $s$, it follows that $f$ is  
linear on the interior of each interval $h(J_i)$,  with slope $\pm s$. Thus
$\tau$ is conjugate to the uniformly piecewise linear map $f$. By the definition of $h$, the conjugacy takes the measure $m$
to Lebesgue measure. 

The converse follows from Example
\ref{1.58}.}

The following is well known for continuous piecewise linear maps; in that case it
follows at once from \cite{Mis-Sz}. See the discussion preceding Proposition \ref{1.9} for the
definition of topological entropy for piecewise monotonic maps that are not continuous.

\prop{1.60}{If $\tau$ is a piecewise linear map with slopes $\pm s$, and $s \ge 1$, then $h_\tau = \ln
s$.}

\prooff{Let $\sigma:X\to X$ be the  local homeomorphism associated with $\tau$, and $\L = \L_\sigma$. Let $R
=\lim_n
\|\L^n 1\|_\infty^{1/n}$, where $\| \,\cdot\, \|_\infty$ denotes the supremum norm.  By \cite[Thm. 6.1]{RuBook},
$R =
 \exp(h_\sigma)$, and by Proposition \ref{1.9}, $h_\sigma = h_\tau$, so to prove the proposition,
we need to show that
$s = R$.

Let $m$ denote Lebesgue measure. Then $m$ is scaled by $\tau$ by the factor $s$. By Proposition \ref{1.56} there is a unique
non-atomic probability measure $\mu$ on
$X$, scaled by $\sigma$ by the factor $s$,  such that
$\mu =
m\circ
\pi$, where
$\pi:X\to I$ is the collapse map. We have
$s^n=\mu(\L^n 1)
\le
\|\L^n 1\|_\infty$, so $s \le \|L^n 1\|_\infty^{1/n}$.
Taking limits, we get $s \le R$.  
Suppose  that $s < R$. Since  $1\le s < R$,  by \cite[Thm. 4.1]{RuBook}, there exists a 
function $\phi$ of bounded variation which is an eigenvector for $\L$ for the
eigenvalue $R$. If $s > 1$, having such an eigenvector   contradicts the spectral radius of $\L$ on the space
$\bv$ of functions of bounded variation being
$s$, cf. \cite[Thm. 1]{Ry}, so $s = R$. (The space $\bv$ is discussed in more detail in the appendix, cf. Definition
\ref{1.113}.) If
$ s = 1$, the entropy of $\tau$ is $0$ (cf.  \cite{BuzziIsom}), so again $s = \exp h_\tau$.}

\cor{1.61}{Let $\tau:I\to I$ be \pwm. If a non-atomic measure $m$ with full support is scaled by $\tau$ by a
factor $s\ge 1$, then $s = \exp(h_\tau)$.}

\prooff{By Proposition \ref{1.59}, $\tau$ is conjugate to a piecewise linear map with slopes $\pm s$. By Proposition
\ref{1.60}, such a map has entropy $\ln s$.}

\section{Decomposition of transitive \pwm\ maps}

In this section we will see that transitive \pwm\ maps often
can  be decomposed into topologically exact
pieces.

In \cite[Thm. 2.5]{Preston}, it is shown
that a transitive, \pwm, continuous map
$\tau$ is either topologically exact, or there exists $c \in I$ such that $\tau$ exchanges $[0,c]$ and
$[c,1]$,  such that  $\tau^2$ is topologically exact on each of $[0,c]$ and $[c,1]$.  We will establish a similar
 finite decomposition  for
 transitive
\pwm\ maps that are not necessarily continuous. To achieve this decomposition, we need to exclude maps
that are ``essentially injective".

\defi{1.62}{A map $\tau:I\to I$ is {\it essentially  injective\/} if there  are
no disjoint intervals $J, J'$ with $\tau(J) = \tau(J')$.}

This is equivalent to $\tau$ being injective on the complement of the set of  endpoints of intervals of
monotonicity.

\lem{1.63}{A \pwm\ map $\tau:I\to I$ is essentially injective iff the associated
local homeomorphism $\sigma:X\to X$ is injective.}

\prooff{\cite[Lemma 11.3]{Shu}}

\prop{1.64}{If $\tau:I\to I$ is \pwm\ and transitive, and
 $\sigma:X\to X$ is the associated local homeomorphism, then there exists a 
scaling measure for $\tau$ and for $\sigma$. Any such scaling measure has  full support and is non-atomic, and
the scaling factor
$s = \exp h_\tau
\ge 1$. Here
$s > 1$ iff
$\tau$ is not essentially injective.}

\prooff{Note that $\sigma$ is strongly transitive (Proposition \ref{1.3}). Let $K$ be the set of probability measures on $X$, viewed
as positive functionals on
$C(X)$. Then
$K$ is a w*-compact convex set.  Define
$T:K\to K$ by
$T\nu = (\nu(\L1))^{-1}\nu\circ
\L$. (Since $\sigma$ is transitive, it is surjective, so
$\L1
\ge 1$, and thus for all $\nu \in K$, 
$\nu(L1) \not= 0$). Note that $T\nu$ is again a probability measure, and
$T$ is a weak* continuous map.  By the Schauder-Tychonov fixed point theorem (cf., e..g, \cite[Thm. V.10.5]{DS}), $T$
has a fixed point in $K$, say $T\mu =
\mu$. Then
$\mu = (\mu(\L 1))^{-1}\mu\circ \L$ so $\mu\circ \L = s\mu$, with $s =
\mu(\L1).$ 

Now let $\mu$ be any probability measure scaled by $\sigma$ by a factor $s$.
Since $\L1\ge 1$, then $s \ge 1$.

Recall that $\sigma$ is injective iff $\tau$ is essentially injective (Lemma \ref{1.63}).  For such maps, $\sigma$ is
invertible, and
$\L_{\sigma^{-1}} = (\L_\sigma)^{-1}$, so $\mu\circ \L_{\sigma^{-1}} = s^{-1}\mu$. Since $\sigma^{-1}$ is
surjective, $\L_{\sigma^{-1}}1\ge 1$, so $s^{-1} \ge 1$, which implies $s =1 $.

On the other hand, if  $\sigma$ is not injective on $X$, suppose that $\sigma(y_1)
= \sigma(y_2) = x$. Then $(\L1)(x) \ge 2$, so $\L1 > 1$ on some open set. Thus $\L1
- 1$ is continuous, nonnegative, and strictly positive on an open set.  We will see in the next paragraph that
$\mu$ has full support, so
$0 < \mu(L1 - 1) = s-1$. Thus $s > 1$.

Now we show that $\mu$ has full support. Let $V$ be any open subset of $X$, and suppose that $\mu(V) = 0$.
By Proposition \ref{1.2.1}, there is a partition of $X$ into clopen
sets on which $\sigma$ is injective. Thus we can write
$V$ as a finite union of open sets on which $\sigma$ is injective, say $V = \cup V_i$.
Then
$\mu(V_i) = 0$ for each $i$, so
\begin{equation}
\mu(\sigma(V)) =
\mu(\sigma(\cup V_i) = \mu(\cup \sigma(V_i)) \le \sum \mu(\sigma(V_i)) = \sum s
\mu(V_i) = 0.
\end{equation}
By the same argument, $\mu(\sigma^k(V)) = 0$ for all $k\ge 0$. Since $\sigma$ is
strongly transitive, a finite number of iterates of $V$ cover $X$, so we conclude that
$\mu(X) = 0$, contrary to $\mu$ being a probability measure. Thus $\supp \mu = X$.

Finally, we will show $\mu$ is non-atomic. If $s > 1$,
then $\mu(\sigma^n(\{x\})) = s^n \mu(\sigma(\{x\}))
< \mu(\sigma(X))$ for all $n$ implies $\mu(\sigma(\{x\})) = 0$, so $\mu$ is non-atomic.  If $s = 1$, then as shown
above, $\sigma$ is injective on
$X$, so it is a homeomorphism. Since $\sigma$ is strongly transitive, it can have no periodic points.  (The
complement of a periodic orbit would be an invariant open set.) For any $x \in X$, $x, \sigma(x),
\ldots$ will be a sequence of distinct points with the same measure, so all must have measure
zero. Thus $\mu$ is non-atomic.}

We will see in Corollary \ref{1.67}  that when $\tau$ is transitive and not essentially injective, then the scaling
measure is unique. 

For any surjective \pwm\ map $\tau:I\to I$, with associated local homeomorphism $\sigma:X\to
X$, the proof above shows that $\tau$ and $\sigma$ possess scaling measures, and shows that for homeomorphisms
the scaling factor must be 1. Thus when $\tau$ is surjective and essentially injective (or equivalently,
$\sigma$ is a homeomorphism), scaling measures are the same as invariant measures.

If $\tau$ is strongly transitive, the following corollary is in \cite{Parry},  with an extra
aperiodicity assumption if $\tau$ is essentially injective.

\cor{1.65}{(Parry)  Let $\tau:I\to I$ be piecewise monotonic and transitive.  Then
$\tau$ is conjugate to a piecewise linear function on
$I$ with slope
$\pm s$, where $s = \exp(h_\tau)$. Here $s = 1$ iff $\tau$ is essentially injective.}

\prooff{By Proposition \ref{1.64}, there exists a non-atomic
scaling measure for  $\tau$, with full support, and with scaling factor $s\ge 1$. The existence of the desired
conjugacy follows from Proposition \ref{1.59}.}

\theo{1.66}{Let $\tau:I\to I$ be \pwm\ and transitive, and not essentially
injective, and let
$\sigma:X\to X$ be the associated local homeomorphism. There is a unique partition of
$X$ into clopen sets
$X_1,
\ldots, X_N$ such that $\sigma(X_i) = X_{i+1 \bmod N}$, with
$\sigma^N|X_i$ topologically exact. Let $\mu$ be a measure scaled by $\sigma$, with scaling factor
$s$.  For each
$i$ there exists
$\phi_i\in C(X)$ of bounded variation,
with
$\supp \phi_i= X_i$ and $\phi_i > 0$ on $X_i$, such that for each $f \in C(X)$ with bounded variation, with $\supp f \subset
X_i$,
\begin{equation}
\lim_k\, (s^{-N}\L_\sigma^N)^kf = \mu(f)\phi_i \hbox{ (uniform convergence)},\label{(1.30)}
\end{equation}
 and such that
the functions $\phi_i$ are cyclically permuted by $(1/s)\L_\sigma$, i.e.,
\begin{equation}
(1/s)\L_\sigma \phi_i =\phi_{i+1\bmod N}.\label{(1.31)}
\end{equation}}

\prooff{By Proposition \ref{1.64}, $s > 1$. We will make use of the results in the appendix. By Corollary \ref{1.121}, there
exists a positive integer
$N$ and disjoint clopen sets
$X_1, \ldots, X_q$ with $\sigma^N$ exact on each $X_i$, and with 
$\sigma(X_i) = X_{\omega(i)}$ for a permutation $\omega$. Since
$\sigma$ is  transitive,  the orbit of each $X_i$ must include
every $X_j$, so 
$\omega$ must be a
cycle, and the  union of all $X_i$ must be $X$. 
By Theorem \ref{1.118} (vi), $q = N$. If necessary, re-index so that $\sigma(X_i) = X_{i+1\bmod N}$. The existence of
$\phi_i \in C(X)$ satisfying  (\ref{(1.30)}) follows from Corollary \ref{1.119} and Proposition \ref{1.120}, and (\ref{(1.31)})
is a consequence of Theorem \ref{1.118} (vi).

To prove uniqueness of this decomposition of $X$, suppose $H_1, \ldots, H_p$ is another partition of $X$ into clopen
sets which are permuted cyclically by $\sigma$,  such that $\sigma^p$ is exact on each $H_i$. Suppose that
$H_i\cap X_j \not= \emptyset$.  Then for $k$ sufficiently large, by exactness of $\sigma^p$ on each $H_i$
and $\sigma^N$ on each $X_j$, $\sigma^{kpN}(H_i\cap X_j) = H_i = X_j$. Thus each $H_i$ coincides with some
$X_j$, and uniqueness follows.}

In the case when $\sigma$ is  positively expansive, the decomposition in Theorem \ref{1.66} (into mixing rather than exact pieces)
follows from
\cite[Thm. 3.4.4]{Aoki-H}.  However, the map $\sigma$ cannot be positively expansive
unless the forward orbit under $\tau$ of the set $C$ of endpoints of intervals of monotonicity
is finite.  Indeed, suppose $x \in \htau C\setminus C$ has an infinite  orbit which never re-enters $C$. Then the points
$\tau^n(x)^\pm$ are all distinct for $n \ge 0$. Viewing $X$ as a compact subset of $\R$,  the sum of the gaps $|\tau^n(x)^+ -
\tau^n(x)^-|$ is finite, so $\lim_n |\sigma^n(x^+) - \sigma^n(x^-)| = \lim_n |\tau^n(x)^+ -
\tau^n(x)^-| = 0$, which shows that $\sigma$ is not positively expansive.

\cor{1.67}{If a \pwm\ map $\tau:I\to I$ is transitive, and is not essentially injective, then there is a
unique scaling measure for
$\tau$, and for the associated local homeomorphism $\sigma:X\to X$.}

\prooff{Let $\mu$ and $\nu$ be scaling measures for $\sigma$. By Proposition \ref{1.64}, both $\mu$ and $\nu$ are
non-atomic and have full support, and have scaling factors $> 1$. Let $X_1, X_2, \ldots, X_N$ be the unique
decomposition given in Theorem \ref{1.66}.  Fix an index $i$, and let $\phi_i$ satisfy (\ref{(1.31)}), where $s$ is the
scaling factor for $\mu$. By (\ref{(1.31)}), $\L_\sigma^N\phi_i=
s^N\phi_i$. Applying $\nu$, $\nu(\L_\sigma^N\phi_i) = s^N \nu(\phi_i)$. Since $\nu$ has full support, then
$\nu(\phi_i) \not= 0$. It follows that the scaling factor for $\nu$ is the same as that for $\mu$.  (This also
follows from Proposition \ref{1.56} and Corollary \ref{1.61}.)

Now from (\ref{(1.30)}), for each index $i$, $\nu(\phi_i) = \mu(\phi_i) = 1$, and if  $f\in C(X)$ has bounded variation
and has  support in 
$X_i$, then
$\nu(f) =
\mu(f)$. Since each clopen subset of $X$ is a finite union of order intervals (Proposition \ref{1.2}), the
characteristic function of each $X_i$ has  bounded variation. If $f\in C(X)$ has bounded variation, then
$f\chi_{X_i}$ is continuous, with bounded variation, and has support in $X_i$. It follows that $\mu$ and $\nu$
agree on all $f \in C(X)$ of bounded variation.  Such functions form a subalgebra of $C(X)$ that is dense by
the Stone-Weierstrass theorem, so $\mu = \nu$.}

If $\tau$ is essentially injective, the uniqueness result above can fail.  For example, there are
minimal interval exchange maps with more than one invariant measure, cf.
\cite{KeaneTwo} or  \cite{Keynes-Newton}.  

The following result also follows from Hofbauer's spectral decomposition for the non-wandering set of \pwm\ maps, 
cf. \cite[Thm. 2]{HofDecomp}, with ``topologically mixing" in place of ``topologically exact".

\cor{1.68}{Let $\tau:I\to I$ be \pwm, transitive, and not essentially injective.  Then
there exist sets
$K_1,K_2
\ldots K_N$ such that
\begin{enumerate}
\item  Each $K_i$ is a finite union of closed intervals.
\item  $\cup_{i=1}^N  K_i= [0,1]$.
\item If $i \not= j$,  the interiors of $K_i$ and $K_j$ are disjoint.
\item For each $i$, $\tau$ maps the interior of $K_i$ onto the interior of
$K_{i+1\bmod N}$
\item For each $i$, $\tau^N$ is topologically exact
when restricted to $K_i$.
\end{enumerate}
\noindent The sets $K_1, \ldots, K_N$ are uniquely determined by these properties.}

\prooff{Let $\sigma:X\to X$ be the associated local homeomorphism, and $\pi:X\to I$ the collapse map. Let
$X_1, \ldots, X_N$ be as in Theorem \ref{1.66}. Define
$K_i = \pi(X_i)$ for $1\le i \le N$. Then (i) follows from the fact that each clopen subset of $X$ is a finite union of order
intervals, and that  $\pi$ maps closed order intervals to closed intervals. Since
$\sigma^N$ is exact on
$X_1,
\ldots, X_n$, then the conjugacy of $\sigma|_{X_0}$ with $\tau|_{I_0}$ implies that $\tau^N$ is exact on each $K_i$.
(See \cite[proof of Lemma 5.2]{Shu} for proof of a similar result.)}

Observe that the decompositions in Theorem \ref{1.66} and
Corollary \ref{1.68} may not exist if $\tau$ is essentially injective. For example,  the map $x
\mapsto x+\alpha \bmod 1$, where $\alpha$ is irrational, is transitive, but is not
topologically exact. There also is no decomposition of some iterate $\tau^N$ into exact pieces, since 
$\tau^N:x\mapsto x+n\alpha
\bmod 1$ is transitive, and thus does not leave invariant any proper closed subset of
$I$.

\defi{1.69}{Let $X$ be a topological space, and $\sigma:X\to X$ a continuous map.  Then
$\sigma$ is {\it topologically mixing\/} if for every pair $U, V$ of non-empty open sets,
there exists $N$ such that $\sigma^n(U) \cap V\not= \emptyset$ for $n \ge N$. (If $\tau:I\to
I$ is \pwm\ but discontinuous at some endpoints, we view $\tau$ as undefined at the set $C$ 
of endpoints of intervals of monotonicity.)}

As in Proposition \ref{1.3} or \cite[Lemma 5.2]{Shu}, one can show that  a \pwm\ map $\tau:I\to I$ will be
topologically mixing iff the associated local homeomorphism $\sigma:X\to X$ has this
property.

Every topologically exact map is topologically mixing. The converse is false, e.g., if $\sigma$ is the
two-sided shift on the space
$\Sigma_2$ of bi-infinite sequences of two symbols, $\sigma$ is topologically mixing, but is not topologically exact. (The
complement of a fixed point is open and invariant.)  On the other hand, for continuous piecewise monotonic maps, being
topologically mixing and being topologically exact are equivalent, cf. \cite[Thm. 2.5]{Preston}.  Furthermore, for
\pwm\ maps that are not essentially injective, mixing and exactness coincide, as we now show.

\prop{1.70}{If $\tau:I\to I$ is \pwm, topologically mixing, and not essentially
injective, then $\tau$ is topologically exact.}

\prooff{Let $\sigma:X\to X$ be the associated local homeomorphism. As observed above, $\sigma$ will be topologically mixing. 
Then
$\sigma$ is transitive, so we can partition $X$ into clopen subsets on which some $\sigma^N$ is exact (Theorem \ref{1.66}).
But if there is more than one member to the partition, this contradicts $\sigma$ being
topologically mixing. So $\sigma$ is exact, and therefore so is $\tau$.}

\section{Simple dimension groups with unique states}

In this section we will see that scaling measures can be used to describe the order on $DG(\tau)$ for
transitive maps $\tau$. In particular, scaling measures induce states (defined below), and states in turn
determine the order.

\defi{1.71}{If $G$ is a dimension group with distinguished order unit $e$, then a {\it state\/}
on
$G$ is a positive homomorphism $\omega$ from $G$ into $\R$ such that $\omega(e) = 1$.  The
{\it strict order\/} on $G$ induced by a state $\omega$ is given by $x \le y$ if either
$\omega(x) < \omega(y)$ or $x = y$.} 

  If $\omega$ is a state on a dimension group $G$, and $\Phi$ is a positive
endomorphism of $G$, we say $\omega$ is {\it scaled by $\Phi$ by the factor $s$} if
$\omega(\Phi(x)) = s\,\omega(x)$ for all $x \in G$.

\prop{1.72}{Let $\tau:I\to I$ be \pwm, and let
$\sigma:X\to X$ be the associated local homeomorphism. If $\mu$ is a probability measure on
$X$ scaled by $\sigma$ by a factor $s$, then the map $[f] \to \mu(f)$ is a state on
$DG(\tau)$. We also denote this state  by $\mu$, and $\mu$ is scaled by $(\L_\sigma)_*$ by the factor $s$.  If $X$ is totally
disconnected, (e.g., if $\tau$ is transitive), this gives a 1-1 correspondence of scaled measures and scaled states.}

\prooff{Clearly $\mu$ is a positive homomorphism from
$C(X,\Z)$ into
$\R$ with value
$1$ on the function constantly $1$. To show $\mu$ is well defined on $DG(\tau)$, suppose $f, g \in C(X,\Z)$
with $[f] = [g]$. Write $\L$ in place of $ \L_\sigma$. Then
$\L^n f = \L^n g$ for some $n \ge 0$, so
\begin{equation}
\mu(f) = s^{-n }\mu(\L^n f) = s^{-n} \mu(\L^n
g) = \mu(g).
\end{equation}
Furthermore, if $[f] \ge 0$, then $\L^n f \ge 0$ for some $n$, so $\mu([f]) = s^{-n}\mu(\L^n f) \ge 0$. 
 Thus $[f]\to \mu(f)$ is a well-defined state on $DG(\tau)$.
Furthermore,
$\mu(\L_* [g]) = \mu(\L g) = s\,\mu(g) = s\, \mu([g])$, so $\mu$ is scaled by $\L_*$ by the factor $s$.

The map that takes a scaled measure to a scaled state is evidently 1-1. Now assume that $X$ is totally disconnected, and
let
$\omega$ be a state on
$DG(\tau)$ scaled by
$\L_*$ by a factor
$s$.  We
will show that $\omega$ comes from a scaled measure.  

Define $\phi:C(X,\Z) \to \R$ by $\phi(f) = \omega([f])$. Then $\phi$
is a positive homomorphism from
$C(X,\Z)$ into
$R$. View $C(X,\Z)$ as a subset of $C(X)$, and let
$C_\Q(X)$ denote the  rational linear span of $C(X,\Z)$ in $C(X)$. Then $\phi$ extends  uniquely
to a positive $\Q$-linear map from $C_\Q(X)$ into $\R$.  Since $\phi(1) = 1$, by positivity
$\phi$ will satisfy
$|\phi(f)|
\le
\|f\|_\infty$ for $f \in C_\Q(X)$. Since
$X$ is totally disconnected, by the Stone-Weierstrass theorem, the (real) linear span of characteristic functions of clopen
subsets of $X$ will be dense in $C(X)$.  Each such real linear combination of characteristic functions can be uniformly
approximated by functions in
$C_\Q(X)$.  It follows that $\phi$ extends uniquely to a positive linear functional on $C(X)$, also denoted by $\phi$. 

Let $\mu$ be the regular
Borel measure such that $\phi(f) = \mu(f)$ for all $f \in C(X)$. Since $\omega$ is scaled by $\L_*$ by the factor $s$, it follows
that
$\phi:C(X)\to \R$ is scaled by $\L$ by the factor $s$, and thus by Lemma \ref{1.55}, $\mu$ is scaled by $\sigma$ by the factor
$s$. Clearly $\omega$ is the scaled state associated with $\mu$, which finishes the proof of the 1-1 correspondence of
scaled measures and scaled states.  Finally, if $\tau$ is transitive, then $X$ is totally disconnected by 
\cite[Prop. 5.8]{Shu},
which completes the proof of the proposition.}

If $\tau:I\to I$ is \pwm, with associated local homeomorphism $\sigma:X\to X$, and $m$ is a non-atomic measure on $I$ with full
support, and $\mu$ is the corresponding measure on $X$ (cf. Proposition \ref{1.56}), then $\mu(I(a,b)) = m([a,b])$ for all $a,
b \in I_1$. Thus we will also refer to the state described in Proposition \ref{1.72} as the state associated with $m$.

If $G, H$ are dimension groups, $\psi:G\to H$ is a {\it dimension group isomorphism\/} if
$\psi $ is both an order isomorphism and a group isomorphism.

\theo{1.73}{Let $\tau:I\to I$ be \pwm\ and transitive, and not essentially injective, with
associated local homeomorphism
$\sigma:X\to X$. Let
$X_1,
\ldots, X_N$ be the partition of $X$ described in Theorem \ref{1.66}, and let
$\mu$ be the unique measure on $X$ scaled by $\sigma$. Then $[f]\mapsto \oplus_i
[f \chi_{X_i}]$ is a unital dimension group isomorphism of $DG(\tau)$ onto $\oplus_i G_i$, where $G_i= \{[f]\in DG(\tau)\mid
\supp f \subset X_i\}$, with distinguished order unit $[\chi_{X_i}]$.  Each $G_i$ is isomorphic to the
dimension group associated with $\sigma^N|X_i$ (cf. Definition \ref{1.3.1}), and is a simple dimension group, with a  unique state
$\mu_i$ given by
$\mu_i([f]) = |\mu(X_i)|^{-1}\mu(f)$, and with the strict order given by that state. Furthermore,
the state on
$DG(\tau)$ induced by $\mu$ is the unique state scaled by
$
(\L_\sigma)_*$, and
$(\L_\sigma)_*$ maps
$G_i$ onto $G_{i+1\bmod N}$ for each $i$. The scaling factor is $s = \exp h_\tau$.}

\prooff{Write $\L$ for $\L_\sigma$.  Since $X_1, \ldots, X_N$ is a partition of $X$ into clopen
$\sigma^N$-invariant subsets,
$(\sigma^N)^{-1}(X_i) = X_i$ for each $i$. It follows that $[f] \mapsto [f|X_i]$ is a group order
isomorphism from $G_i$ onto $G_{\sigma^N|X_i}$.  Since $\sigma$ is a
piecewise homeomorphism, then  $\sigma^N$ is also a piecewise homeomorphism, so each $G_i$  is simple by
exactness of $\sigma^N|X_i$, cf. \cite[Theorem 5.3]{Shu}. 
 From this the fact that $[f]\mapsto \oplus_i [f|X_i]$ is a
dimension group isomorphism is immediate.

Fix an index $i$, with $1 \le i \le N$. We will show that $G_i$ has a 
unique state given by $\mu|X_i$, and has the strict order given by that state. 
By  Theorem \ref{1.66},
there is a continuous function $\phi_i$, strictly positive on $X_i$, such that for each $f \in C(X)\cap \bv$
with support in
$X_i$,
\begin{equation}
\lim_{k\to \infty} s^{-Nk}L^{Nk} f = \mu(f)\phi_i,
\end{equation}
with uniform convergence. Note that each function in $C(X,\Z)$ has bounded variation (since each
clopen set in $X$ is a finite union of order intervals, cf. Proposition \ref{1.2}). Now suppose
$f, g
\in C(X,\Z)$ have support in
$X_i$, and
$\mu(f) < \mu(g)$. Choose $\epsilon> 0$ so that $\epsilon \le (\mu(g) - \mu(f))\phi_i(x)$ for
all
$x \in X_i$, so that
$\mu(f)\phi_i \le \mu(g)\phi_i -
\epsilon$. Choose
$n$ so that 
\begin{equation}
\| s^{-n}\L^n f - \mu(f)\phi_i\|_\infty < \frac{\epsilon}{2}
\text{ and }  \| s^{-n}\L^n g - \mu(g)\phi_i\|_\infty< \frac{\epsilon}{2}.
\end{equation}
 Then 
\begin{equation}
s^{-n}\L^n f \le \mu(f)\phi_i +
\frac{\epsilon}{2} \le \mu(g)\phi_i - \frac{\epsilon}{2} \le s^{-n}\L^n g.
\end{equation}
Thus $\L^n f \le \L^n g$, which implies $[f] \le [g]$. 

Conversely, suppose $[f] \le [g]$. Choose $n$ so that $\L^nf \le
\L^n g$. Then 
\begin{equation}
\mu(f) = \mu(s^{-n}\L^n f) \le \mu(s^{-n}\L^ng)= \mu(g).
\end{equation}
If $\mu(f) = \mu(g)$, then $\mu$ is zero on the nonnegative continuous function
$\L^ng-\L^nf$. Since $\supp \mu = X$, then $\L^n f = \L^n g$, so $[f] = [g]$. Thus
either $\mu(f) < \mu(g)$ or else $[f]=[g]$, which completes the proof that the order on
$G_i$ is the strict order induced by the state $\mu_i$.

To show there is a unique state on $G_i$, suppose $\omega$ is a state
on $G_i$, and fix $f \in C(X,\Z)$ with support in $X_i$. Let $\epsilon > 0$ be given, and choose a
positive integer $n$ such that $2/n < \epsilon$. Now chose integers $k$ and $p$ such
that $|\mu_i(f) - (k/p)| < 1/n$. This is equivalent to
\begin{equation}
nk-p < np\,  \mu_i(f) < nk + p.
\end{equation}
Since the ordering on $G_i$ is the strict ordering given by $\mu_i$, it follows that
\begin{equation}
(nk-p)[\chi_{X_i}] < np\, [f] < (nk + p)[\chi_{X_i}].
\end{equation}
Applying the state $\omega$ gives
\begin{equation}
nk - p < np\, \omega([f]) < nk+ p.
\end{equation}
Dividing by $np$ gives
\begin{equation}
|\omega([f])-  (k/p)| \le \frac{1}{n}.
\end{equation}
Thus 
\begin{equation}
|\omega([f]) - \mu_i(f)| \le 2/n< \epsilon.
\end{equation}
Since $\epsilon > 0$ was arbitrary, we conclude that $\omega([f]) = \mu_i(f)$, and thus
$\omega $ is the state induced by $\mu_i$, proving that $G_i$ has a unique state.

Since $\sigma(X_i) = X_{i+1\bmod N}$, then
\begin{equation}
\supp f \subset X_i \implies \supp (\L f) \subset \sigma(\supp [f]) \subset X_{i+1\bmod N}.
\end{equation}
It follows that
$\L_*G_i
\subset G_{i+1\bmod N}$. Since $\sigma$ is surjective, then $\L_*$ is
surjective on $G$ (\cite[Lemma 3.5]{Shu}), so we must have $\L_*G_i
= G_{i+1\bmod N}$ for each $i$.

Finally, suppose that $\mu$ is scaled by the factor $s$, and let $\omega$ be any state on $DG(\tau)$ that
is scaled by a factor
$s_1$. Fix an index $i$. By the uniqueness of states on $G_i$, there is a constant $k$ such that $\omega
= k\mu$ on
$G_i$. Then
$\omega
\circ \L_*^N =  k \mu\circ \L_*^N$ on $G_i$, so $s_1^N \omega = k s^N \mu$ on $G_i$. If $\omega|G_i$ were
0, then the fact that $\omega$ is scaled by $\L_*$ and that $G_1, \ldots, G_N$ are cyclically permuted by
$\L_*$ would imply that $\omega = 0$.  Therefore we conclude that $s_1 = s$. Now it follows that
$\omega\circ \L_*^j = k \mu\circ \L_*^j$ on $G_i$ for $j = 1, \ldots, N$, so $\omega = k\mu$ on all $G_j$. Then $1 = \omega(1) =
k\mu(1)$ implies that 
$k = 1$, so $\omega = \mu$.

The fact that the scaling factor is $s = \exp h_\tau$ follows from Corollary \ref{1.61}.}

Theorem \ref{1.73} provides a description of all states on $DG(\tau)$.  For each $i$, define a state
$\widetilde \mu_i$ on $DG(\tau)$ by $\widetilde \mu_i([f]) = \mu(X_i)^{-1}\mu(f\chi_{X_i})$. Then
$\widetilde \mu_i| G_i = \mu_i$, and every state on
$DG(\tau)$ is a unique convex combination of
$\widetilde\mu_1,
\ldots,
\widetilde\mu_N$, so the state space is a simplex with $\widetilde\mu_1, \ldots, \widetilde\mu_N$ as extreme
points.

\cor{1.74}{If maps $\tau_i:I\to I$ are \pwm\ and transitive for $i = 1,2$, and if the dimension triples for
$\tau_1$ and $\tau_2$ are isomorphic, then $\tau_1$ and $\tau_2$ have the same topological entropy.}

\prooff{Assume that the dimension triples for $\tau_1$ and $\tau_2$ are isomorphic. Suppose first that neither map is
essentially injective. Then each dimension group has a unique scaled state (Theorem \ref{1.73}). By the isomorphism of the
dimension triples, the scaling factors must coincide. Since the scaling factors are given by the exponential of the
topological entropy, the two entropies must coincide.

Suppose now that $\tau_1$ is essentially injective, and that $\sigma_1:X^1\to X^1$ is the associated local homeomorphism.
By Proposition \ref{1.64},
$\tau_1$ admits a scaling measure, and all such scaling measures have scaling factor $s = \exp h_{\tau_1} = 1$, so
$h_{\tau_1} = 0$. By Proposition \ref{1.72}, all scaled states on $DG(\tau_1)$ have scaling factors 1. The same must hold for
$DG(\tau_2)$. By Proposition \ref{1.64} and Proposition \ref{1.72}, this implies that $\tau_2$ is essentially injective, and so has
zero topological entropy. In particular, $h_{\tau_1} = h_{\tau_2}$.}

\cor{1.75}{If $\tau:I\to I$ is \pwm\ and topologically exact, then $DG(\tau)$ is simple and has a unique state, given by
the unique scaling measure, and has the strict order given by that state.}

\prooff{If $\tau$ is topologically exact, then so is the associated local homeomorphism $\sigma:X\to X$ (Proposition \ref{1.3}).  It
follows that $\sigma$ is not injective, and so $\tau$ is not essentially injective (Lemma \ref{1.63}). Now the
conclusion follows from Theorem
\ref{1.73}.}

Renault \cite[Example 6.2]{RenRN} has proven a similar uniqueness result for states on dimension groups associated with positively
expansive maps that are topologically exact. (See the remarks on positively expansive maps
following Theorem
\ref{1.66}.)

\defi{1.76}{Let $G$ be a dimension group with distinguished order unit $e$. An element $x\in G$ is {\it
infinitesimal\/} if for all
$n \ge 0$,
$- e
\le n [f] \le e$.  (This is independent of the choice of order unit.) The subgroup of infinitesimals
is denoted
$G_{\rm inf}$.}

Observe that if $G$ is a dimension group with distinguished order unit $e$, and $\Phi:G\to G$ is a
positive endomorphism, then
$G_{\rm inf}$ will be invariant under $\Phi$.

If $G$ is an ordered abelian group, and $H$ is a subgroup that is an order ideal (i.e., $g \in G$ and $ 0
\le g
\le h
\in H$ imply $g \in H$), then we define  $(G/H)^+$ to be the image under the quotient map of
the positive cone of $G$. With the  ordering given by $(G/H)^+$, $G/H$ becomes an ordered
abelian group.

Below we view
${\R}^N$ as an ordered group with the usual coordinatewise order; it is evidently a dimension group.

\prop{1.77}{Let $\tau:I\to I$ be \pwm\ and transitive and not essentially injective, let $\sigma:X\to X$ be
the associated local homeomorphism, and let
$DG(\tau) = \oplus_i G_i$ be the decomposition given in Theorem \ref{1.73}. Let
$\mu$ be the unique measure on $X$ scaled by $\sigma$, with scaling factor $s$.  Then we have an
isomorphism of dimension groups given by
\begin{equation}
DG(\tau)/DG(\tau)_{\rm inf} \cong \oplus_i G_i/(G_i)_{\rm inf} \cong \oplus_i \mu(G_i),\label{(1.32)}
\end{equation}
where each $G_i/(G_i)_{\rm inf}$ has the quotient ordering, and $\oplus_i \mu(G_i)$ has the ordering
inherited from $\R^N$.  An isomorphism from $DG(\tau)/DG(\tau)_{\rm inf}$ onto $\oplus_i \mu(G_i)$ is given by
$\Psi([f]) =\oplus_i \mu(f\chi_{X_i})$, and $\Psi$ carries the order automorphism induced by $\L_*$ on
$DG(\tau)/DG(\tau)_{\rm inf}$ to the order automorphism of $\oplus_i \mu(G_i)$ given by
\begin{equation}
(\lambda_1, \lambda_2, \ldots
\lambda_N) \mapsto (s\lambda_N, s\lambda_1, \ldots, s\lambda_{N-1}).
\end{equation}} 

\prooff{If $X_1, \ldots, X_N$ is the  partition of $X$ given in Theorem \ref{1.66},
then $[f]
\in DG(\tau)$ is an infinitesimal in $DG(\tau)$ iff each $[f \chi_{X_i}]$ is infinitesimal in $G_i$, so
$G_{\rm inf} = \oplus_i (G_i)_{\rm inf}$. Since $G = \oplus_i G_i$,  the first isomorphism in (\ref{(1.32)}) follows.
For the second isomorphism, since
$G_i$ has the strict ordering from
$\mu|G_i$, then
$[f]  \mapsto \mu(f)$ is a positive homomorphism from $G_i$ onto $\mu(G_i)$ with kernel $(G_i)_{\rm inf}$,
and the quotient ordering is exactly the ordering on $\mu(G_i)$ inherited from $\R$. 

 To verify the final statement of the proposition, for $f \in C(X,\Z)$ and $1\le i \le N$, let $f_i =
f\chi_{X_i}$, and $\lambda_i = \mu(f_i)$. Then $\Psi([f]) = (\lambda_1, \ldots, \lambda_N)$.
Since $\supp \L f_i \subset X_{i+1\bmod N}$ for all $i$, then the $k$th coordinate of $[\L f]$ in $\oplus_i G_i$ is $[\L f_{k-1\bmod
N}]$, so
\begin{equation}
(\Psi([\L f]))_k  = \mu(\L f_{k-1\bmod N}) = s \mu(f_{k-1\bmod N}) = s\lambda_{k-1\bmod N},
\end{equation}
which completes the proof of the proposition.}

Krieger proved a  result analogous to Proposition \ref{1.77} with two-sided irreducible shifts of finite type in place of transitive
\pwm\ maps $\tau$, and with the invariant measure of maximal entropy in place of the scaling measure $\mu$ (\cite[Thm.
3.2]{KriShift}). He then used the range of this measure on the clopen subsets  as an invariant to distinguish certain shifts of
finite type. The desire to generalize this result to interval maps was one motivation for the current paper.

\section{Dimension triples as invariants for \pwm\ maps}

The dimension triple  determines the topological
entropy (Corollary \ref{1.74}). However, to determine $\tau$ up to conjugacy, the dimension triple is still not quite
enough, and we will now discuss some additional information that suffices.

\defi{1.78}{Let $\tau:I\to I$ be \pwm, continuous, and transitive, with associated partition $0
= a_0 < a_1 < \ldots < a_n= 1$. The $n$-tuple $(I(a_0,a_1), I(a_1,a_2), \ldots, I(a_{n-1},a_n))$ is
called the {\it canonical sequence of generators for the dimension module $DG(\tau)$}.  We will call $(a_0,a_1)$  the {\it
first interval of monotonicity}.}

\theo{1.79}{Let $\tau_1:I\to I$ and $\tau_2:I\to I$ be continuous, transitive, \pwm\ maps. Then there is an
increasing conjugacy from $\tau_1$ to
$\tau_2$  iff
\begin{enumerate} 
\item the maps $\tau_1$ and $\tau_2$ are both increasing, or both decreasing, on their first interval
of monotonicity, and
\item there is a unital isomorphism $\Phi$ from the dimension triple for $\tau_1$ onto the
dimension triple for  $\tau_2$, taking the canonical sequence of generators for $DG(\tau_1)$ onto that for
$DG(\tau_2)$.
\end{enumerate}}

\prooff{Assume that (i) and (ii) hold. For $i = 1,2$, let 
$\sigma_i:X^i\to X^i$ be the  local homeomorphism associated with $\tau_i$.  Note that since each
$\tau_i:I\to I$ is continuous and transitive, it is surjective, and is not essentially injective. (If it
were essentially injective, by continuity it would be injective, and thus would be a homeomorphism. However,
no homeomorphism of $I$ is transitive.)  Thus by Corollary \ref{1.67} and Proposition \ref{1.56}, for $i = 1,2$, there is a 
unique measure $\mu_i$ on
$[0,1]$ scaled by
$\tau_i$ by a factor $s_i$, with associated scaled measure $\widetilde \mu_i$ on $X^i$.  Define
$h_i:[0,1]
\to [0,1]$ by
$h_i(x) =\mu_i([0,x])$. Then
$h_i$ is an increasing conjugacy from $\tau_i$ onto the uniformly piecewise linear map $T_i = h_i\circ
\tau_i\circ h_i^{-1}$, with slopes $\pm s_i$ for $i = 1,2$. (See the proof of Proposition \ref{1.59}). We will show that
$T_1 = T_2$,  and then
$h_2^{-1}\circ h_1$ will be the desired increasing conjugacy from $\tau_1$ to $\tau_2$.

 If $a_0,
a_1,
\ldots, a_n$ and
$b_0, b_1,
\ldots, b_n$ are the  partitions associated with $\tau_1$ and $\tau_2$ respectively, then the lengths
of the
$i$-th interval of monotonicity for $T_1$ and $T_2$ will be
$\mu_1(a_{i-1},a_i) = \widetilde \mu_1(I(a_{i-1},a_i))$ and $\mu_2(I(b_{i-1},b_i)) =  \widetilde
\mu_2(I(b_{i-1},b_i))$ respectively. By the uniqueness of the scaled states on
$DG(\tau_1)$ and $DG(\tau_2)$, $\Phi$ must take $\widetilde \mu_1$ to $\widetilde \mu_2$, and we must have
$s_1 = s_2$. Since $\Phi$ takes the canonical sequence of generators of $DG(\tau_1)$ to the corresponding sequence for
$DG(\tau_2)$, then  
$\widetilde\mu_1(I(a_{i-1},a_i))=
\widetilde\mu_2(I(b_{i-1},b_i))$ for all $i$. Thus $T_1$ and $T_2$ have intervals of monotonicity of the same
length, increase and decrease on these intervals in the same order, and have the same slopes. Thus for all $x$, we
have $T_2(x) = T_1(x) + T_2(0) - T_1(0)$. Since $\max_{x\in I} T_1(X) = \max_{x\in I} T_2(x) = 1$, we must have $T_1(0) = T_2(0)$.
Thus  $T_1 = T_2$, which completes the proof that $\tau_1$ and $\tau_2$ are conjugate.

Conversely, if there is an increasing conjugacy $h$ from $\tau_1$ to $\tau_2$, then (i) will hold, and
the conjugacy will  carry each interval $(a_{i-1}, a_i)$ to the corresponding interval $(b_{i-1},b_i)$. From
the construction of the associated local homeomorphisms $\sigma_i:X^i \to X^i$, $h$ will lift to a conjugacy
$\widetilde h$ from $(X^1,\sigma_1)$ onto $(X^2,\sigma_2)$, carrying $I(a_{i-1},a_i)$ to $I(b_{i-1},b_i)$.
Then $[f] \mapsto [f\circ \widetilde h^{-1}]$ will be the desired isomorphism of dimension group triples.}

We will see in Example \ref{1.98} that condition (i) in Theorem \ref{1.79} is not redundant.
\medskip

Define $\phi:I\to I$ by  $\phi(x) = 1-x$. If $\tau_1:I\to I$ is any \pwm\ map, and $\tau_2 = \phi\circ
\tau_1\circ \phi^{-1}$, then $\phi$ will induce an isomorphism from $DG(\tau_1)$ onto $DG(\tau_2)$ that will
carry the sequence of canonical generators of $DG(\tau_1)$ onto the sequence of canonical generators of
$DG(\tau_2)$ in reverse order.   Since every
decreasing homeomorphism is the composition of an increasing homeomorphism with $\phi$, it isn't difficult to
reformulate Theorem \ref{1.79}  to include arbitrary (increasing or decreasing) conjugacies; we leave that to the
reader.

For particular families of interval maps, it can happen that the dimension triple 
determines the map
among members of that family.  For example, we will see later that transitive unimodal maps are determined
up to conjugacy by their entropy, which, as we have seen, can be recovered from the dimension triple.

For a related discussion of dimension groups as an invariant for minimal  homeomorphisms of a Cantor set,
see \cite{GPS}, where it is shown, for example, that two such homeomorphisms are orbit equivalent
iff certain associated  dimension groups modulo infinitesimals are unitally order isomorphic. 

\begin{figure}[htb]
\centerline{\includegraphics[scale=0.90]{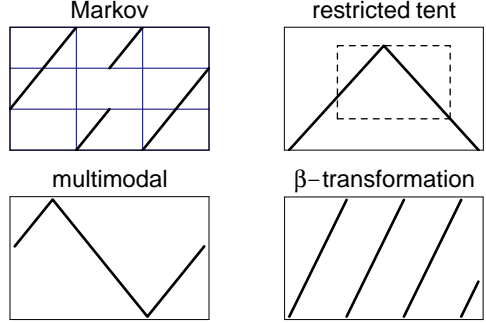}}
\caption{\label{fig1} Examples of transitive piecewise monotonic maps}
\end{figure}

In the remainder of this paper, we will use the  results of previous sections to compute the dimension groups for several families
of transitive maps.  For later reference, some of these are illustrated in Figure \ref{fig1}.

\section{Markov maps}

For use in later sections, we review the
description of
$DG(\tau)$ for Markov maps
$\tau$ from section 8 of
\cite{Shu}. A
\pwm\ map
$\tau:I\to I$ is Markov if there is a partition $B = \{0= b_0, b_1, \ldots, b_n = 1\}$ of $I$ such that $\tau$ is monotonic
on each interval
$[b_{i-1}, b_i]$, and such that the image of each such interval is a union of intervals of the form $[b_{j-1},b_j]$.  We
also require that the orbits under
$\htau$ of the critical points and endpoints eventually land in
$B$.  Observe that any Markov map is necessarily eventually surjective, since the successive images of the whole space decrease
and each image is  a union of some of the Markov partition intervals. 

\defi{1.81}{If $A$ is an $n \times n$ \zeroOne  matrix, $G_A$ denotes the stationary inductive limit
 $\Z^n \mapright{A}\Z^n$ in the category of ordered abelian groups. We represent members of $G_A$ as equivalence classes $[v,n]$
where
$v
\in \Z^n$ and $n \in N$, with $[v_1,n_1] = [v_2,n_2]$ if $v_1A^{n_2+k} = v_2A^{n_1+k}$ for some $k\ge 0$. 
$G_A$ will be a dimension group.  The action of $A$ induces an automorphism of $G_A$ denoted $A_*$, where $A_*[v,n] = [vA,n]$,
so we view
$G_A$ as a
$\Z[t,t^{-1}]$ module, and refer to $(G_A,G_A^+, A_*)$ as the dimension triple associated with $A$.}

\prop{1.82}{(\cite[Prop. 8.4]{Shu}) Let  $\tau:I\to I$ be \pwm\ and Markov, with associated local homeomorphism $\sigma:X\to X$.
Let
$E_1, E_2,
\ldots, E_q$ be the associated Markov partition for $\sigma$, with incidence matrix
$A$, and define $\psi:\Z^q \to C(X,\Z)$ by $\psi(z_1, z_2, \ldots, z_{q}) = \sum_i z_i E_i$.  Then the map 
$\Phi:G_A\to DG(\tau)$ defined by
$\Phi([v,n]) = \L_*^{-n}[\psi(v)]$ is an isomorphism from the dimension triple
$(G_A,G_A^+,A_*)$ onto the dimension triple $(DG(\tau),DG(\tau)^+,\L_*)$.}

\exe{1.87}{Let $\tau$ be the Markov map in Figure \ref{fig1}, and let $\sigma:X\to X$ be the associated local
homeomorphism. Then $\tau$ is Markov, and the associated incidence matrix
$A$ is primitive. If $(X_A,\sigma_A)$ is the associated shift of finite type, by \cite[Cor. 8.9 and Ex. 8.10]{Shu},
$(X,\sigma)$ is conjugate to
$(X_A,
\sigma_A)$, and the dimension triples
$(DG(\tau),DG(\tau)^+, \L_*)$ and $(G_A, G_A^+, A_*)$ are isomorphic. The unique measure scaled by $\tau$ is Lebesgue
measure $m$, with scaling factor $2$. If $E_1 = I(0,1/3)$, $E_2= I(1/3,2/3)$, $E_3= I(2/3,1)$ are the three intervals
of the Markov partition, then the image of the unique state is $\sum_i\Z[s,s^{-1}]m(E_i) = (1/3) \Z[1/2]$. By Proposition
\ref{1.77}, we also have $DG(\tau)/DG(\tau)_{\rm inf} \cong (1/3) \Z[1/2]$.}

When $\tau$ is \pwm, topologically exact, and Markov, with incidence matrix $A$, it is interesting to know whether the unique state 
gives an isomorphism of
$DG(\tau)
\cong G_A$ onto the range of the state, i.e., when there are no infinitesimals.  By \cite[Thm. 5.10 and Cor. 5.11]{BMT}, this
happens iff the characteristic polynomial
$p(t)$ of the primitive matrix
$A$, with the  largest power of $t$ divided out, is irreducible. In Example \ref{1.87}, the characteristic polynomial of the incidence
matrix $A$ is $p(t) = t^3 - 3t-2= (t+1)(t^2-t-2)$. Since $p$ is not irreducible, then there are infinitesimals.

\section{Unimodal maps}

A \emph{unimodal map} $\tau:I\to I$ is a continuous piecewise monotonic map with just two intervals of monotonicity $[0,c]$ and
$[c,1]$, such that $\tau$ increases on $[0,c]$ and decreases on $[c,1]$. In this section we will describe all transitive
unimodal maps, and compute their dimension triples.

Given $s> 1$, we define the {\it restricted tent map\/} $T_s$ by 
\begin{equation}
T_s(x) = \begin{cases}
1+s(x-c) & \text{if $x \le c$}\cr
1-s(x-c) &\text{if $x > c$},
\end{cases}
\end{equation}
where $c = 1-1/s$. 
This is the usual symmetric tent map $\tau_0$ on $[0,1]$ with slopes $\pm s$, restricted to the interval $[\tau_0^2(1/2), \tau_0
(1/2)]$, which is the interval of most interest for the dynamics.    Then the map has been rescaled so that its domain
is  [0,1] (See Figure \ref{fig1}.) Note that $T_s(c) = 1$ and $T_s(1) = 0$.
The  following is essentially in \cite{Jonk-Rand} and \cite[p. 245]{Melo-S}.

\lem{1.89}{Let $ T = T_s$ be the restricted tent map, with slopes $\pm s$. 
\begin{enumerate}
\item If  $\sqrt{2} < s \le 2$, then $T$ is topologically mixing (and exact).
\item If $s = \sqrt{2}$, and $p$ is the fixed point of $T$, then $T$ exchanges
$[0,p]$ and $[p,1]$, and $T^2$ is conjugate to the full tent map on each of these intervals. 
Thus $T$ is transitive, but is not topologically exact.
\item If $s < \sqrt{2}$, then $T$ is not transitive.
\end{enumerate}}

\prooff{Let $c = 1-1/s$, so that $c$ is the critical point of $T = T_s$, and let $p$ be the
fixed point of $T$.  If
$s <
\sqrt{2}$, then
$0 < c < T^2(0) < p < T(0)< T^3(0) < 1$, so that $[0,T^2(0)] \cup [T (0),1]$ is a closed
invariant subset, with non-empty interior, contradicting transitivity. If $\sqrt{2} < s \le 2$, then for any interval
$J$, either $J$ and $T(J)$ contain  $c$ (in which case $T^2(J)$ contains 0 and 1, so $T^2(J) = [0,1]$), or
else applying $T^2$ to $J$ gives an interval whose length is at least $s(s/2) = s^2/2> 1$ 
times the length of $J$.  Thus for some $n$, $T^n(J) = [0,1]$. If $s = \sqrt{2}$, then $T(0)$ is fixed by $T$, and the
statements in the proposition for
$s=
\sqrt{2}$ are straightforward to check.}

 Recall that $\Z[t,t^{-1}]$ denotes the additive group  of Laurent polynomials
  with integer coefficients. The group $\Z[t,t^{-1}]$ has no given order. For each positive $s \in \R$, 
we write $\widehat s:\Z[t,t^{-1}]\to \R$ for the evaluation map, i.e., $\widehat s(p) = p(s)$.  For each such $s$, we  define a (partial) order on $\Z[t,t^{-1}]$ 
by  $p\ge 0$ iff $p = 0$ or $p(s) > 0$, and  call this order the \emph{strict order given by evaluation at $s$}.

 For $0 < s \in \R$, we denote by $\Z[s,s^{-1}]$ the subgroup of $\R$ generated by $\{s^n \mid n \in \Z\}$, or equivalently, the
subring of $\R$ generated by $s$ and $s^{-1}$.
 We equip $\Z[s,s^{-1}]$ with the order inherited from $\R$.

We summarize results from \cite[Thm. 9.1]{Shu} for $DG(\tau)$ when $\tau$ is unimodal and surjective. 
Analysis of dimension groups for unimodal maps divides naturally into the cases where the map is Markov or not.
 If
the critical point
$c$ is eventually periodic, then $\tau$ will be Markov, so $(DG(\tau), DG(\tau)^+, \L_*)$ will be isomorphic to $(G_A,
G_A^+, A_*)$. If instead the orbit of
$c$ is infinite, then $p \mapsto p(\L_*)I(0,1)$ is an isomorphism of  $\Z[t,t^{-1}]$ onto $DG(\tau)$ as an abelian
group, with the action of
$\L_*$ on $DG(\tau)$ corresponding to multiplication by $t$ on $\Z[t,t^{-1}]$. (Here I(0,1) is the class in $DG(\tau)$ of the
function identically 1 on $X$; see the remark before Theorem 2.12.)

If $\tau$ is transitive, we now
describe the order on
$DG(\tau)$ in the case where
$c$ is not eventually periodic.

\prop{1.90}{Let $\tau:I\to I$ be unimodal and transitive, and let $s = \exp h_\tau$.
Then $\sqrt{2} \le s \le 2$, and $\tau$ is conjugate to the restricted tent map $T_s$. 
If in addition the critical point $c$ is not
eventually periodic, then
\begin{enumerate}
\item The map $p\mapsto p(\L_*)I(0,1)$ is a group isomorphism of $ \Z[t,t^{-1}]$ onto  $DG(\tau) $. 
 This isomorphism carries the evaluation map $\widehat s:\Z[t,t^{-1}]\to \R$ to the
unique scaled state on $DG(\tau)$.
\item If $s > \sqrt{2}$, and we give  $\Z[t,t^{-1}]$ the strict order given by
evaluation at
$s$, then the group isomorphism in (i) is also an order isomorphism. 
\item   If $s > \sqrt{2}$ then the unique state on $DG(\tau)$ induces a group and order isomorphism from $DG(\tau)/DG(\tau)_{\rm
inf}$ onto
$\Z[s,s^{-1}]\subset
\R$, with the automorphism of $DG(\tau)/DG(\tau)_{\rm inf}$ induced by $\L_*$ carried to  multiplication by $s$.
\end{enumerate}}

\prooff{By Corollary \ref{1.65}, $\tau$ is conjugate to a piecewise linear map $T$ with slopes
$\pm s$, with $s = \exp h_\tau$.   If
$T(1) > 0$, then $T(0) = 0$, and for $\epsilon > 0$ small enough, and $J = (\epsilon/2,\epsilon)$,
$T^n(J)
\cap J =
\emptyset$ for all $n \ge 1$. This would imply that $\tau$ is not transitive, so we must
have 
$T(1) = 0$. Then $T = T_s$, and $\sqrt{2} \le s \le 2$ by Lemma \ref{1.89}.

Assume now that the orbit of $c$ is not eventually periodic.  As discussed in the remarks preceding this
proposition, the map $p\to p(\L_*)I(0,1)$ is an isomorphism of $\Z[t,t^{-1}]$ onto $DG(\tau)$.  
Since $\tau$ is
transitive, there is a unique scaled state $\omega$ on $DG(\tau)$ (Theorem \ref{1.73}).  By the scaling property,
 $\omega(p(\L_*)I(0,1)) = p(s) = \widehat s(p)$. Thus the isomorphism $p \mapsto p(\L_*)I(0,1)$ carries $\widehat s$ to
$\omega$. 

When $\sqrt{2}< s \le 2$, 
then  $T_s$ is topologically exact (Lemma \ref{1.89}), and since $\tau$ is conjugate to $T_s$, then $\tau$ is topologically exact. 
Hence by Corollary \ref{1.75} there is just one state on $DG(\tau)$, and $DG(\tau)$ is a simple dimension group with the strict 
order from this state.  Therefore by (i), if we give $\Z[t,t^{-1}]$ the strict order given by evaluation at $s$, then $p \mapsto p(\L_*)I(0,1)$ 
becomes an order isomorphism as well as a group isomorphism.

 By exactness of $\tau$ and Proposition \ref{1.77},  the unique state $\omega$ on $DG(\tau)$ is a group and order isomorphism
of 
$DG(\tau)/DG(\tau)_{\rm inf}$ onto $\omega(DG(\tau)) \subset \R$. Since this unique state must coincide with
the unique scaled state, which has scaling factor $s$, for $p \in \Z[t,t^{-1}]$ we have
$\omega(p(\L_*)I(0,1)) = p(s)$, so $\omega(DG(\tau)) \supset \Z[s,s^{-1}]$. 
Since $p\to p(\L_*)I(0,1)$ maps $\Z[t,t^{-1}]$ onto $DG(\tau)$, the image of $\omega$ equals $\Z[s,s^{-1}]$. 
The fact that $\omega$ carries the action of $\L_*$ to multiplication by $s$ is just a restatement of the fact that $\omega$ is
scaled by $\L_*$ by the factor $s$.}

In Proposition \ref{1.90}, if $s = \sqrt{2}$, the dimension group for $T_s$ is the direct sum of two copies of the dimension group
for
$T_2$. (See Lemma \ref{1.89} and Theorem \ref{1.73}).  The dimension triples for $T_2$ and $T_{\sqrt{2}}$ are computed in
\cite[Examples 9.2, 9.3]{Shu}.

\cor{1.91}{Two unimodal transitive maps are conjugate iff their  dimension triples
are isomorphic.}

\prooff{By Corollary \ref{1.74}, two transitive \pwm\ maps with   isomorphic dimension triples have the same topological
entropy. By Proposition  \ref{1.90}, transitive unimodal maps with the same topological entropy are conjugate.}

\exe{1.93}{Let $\tau$ be the restricted tent map $T_s$, with $\sqrt{2} < s < 2$.  Since $\tau$ is topologically exact, then
$DG(\tau)$ is simple (Corollary \ref{1.75}).  If
$s$ is transcendental, then 
$p(s)$ is non-zero for all $p
\in
\Z[t,t^{-1}]$, so the unique state is an order automorphism of
$DG(\tau) = \Z[t,t^{-1}]I(0,1)$ onto
$\Z[s,s^{-1}]\subset \R$.}

\exe{1.94}{Let $s = 3/2$. Then the critical point of the restricted tent map $T_s$  is not
eventually periodic, cf. \cite[Example 9.4]{Shu}. By Proposition \ref{1.90},
$DG(T_s)
\cong
\Z[t,t^{-1}]$, with the strict order given by the unique state, namely evaluation at $s$.  The image of this state is
$\Z[s, s^{-1}]$, but this state is not 1-1, e.g. $p(t) = 2t-3$ is sent to zero. Thus
there are infinitesimals in this dimension group.} 

\section{Multimodal maps}

In \cite[Prop. 10.6]{Shu}, sufficient conditions are given for a continuous $n$-modal map 
to have a dimension group isomorphic as an (unordered) module to $(\Z([t,t^{-1}])^{n-1}$. 
We will see that if we assume the maps are mixing and uniformly piecewise linear, we can give an explicit
description of the (ordered) dimension group modulo the subgroup of infinitesimals.  With some extra assumptions that
guarantee that there are no infinitesimals, we can  describe explicitly the dimension module itself.  Since every
continuous mixing \pwm\ map is conjugate to a uniformly piecewise linear map, these results are also applicable to such
maps. 

Recall that a continuous topologically mixing \pwm\ map is  topologically exact, cf. \cite[Thm. 2.5]{Preston} or
Proposition \ref{1.70}.

\prop{1.96}{Let $\tau:I\to I$ be continuous, piecewise linear, and topologically mixing,
with slopes
$\pm s$,  and with associated partition $a_0 < a_1 < \cdots < a_n$.  Then $DG(\tau)/DG(\tau)_{\rm inf} \cong \sum_{i=1}^n
\Z[s,s^{-1}](a_i-a_{i-1})$ (with the order inherited from $\R$, and with the action of $\L_*$ given by multiplication by
$s$).}

\prooff{By Corollary \ref{1.75}, there is a unique state, which is given by the unique scaling
measure, namely, Lebesgue measure. By Corollary \ref{1.42}, the intervals $I(a_{i-1},a_i)$ for $1\le i\le n$ generate the
dimension module $DG(\tau)$. By Proposition \ref{1.77}, 
the ordered group $DG(\tau)/DG(\tau)_{\rm inf}$ is isomorphic to the range of that state, namely,
$\sum_i\Z[s,s^{-1}](a_i-a_{i-1})$.}

\prop{1.97}{Let $\tau:I\to I$ be continuous, piecewise linear, and topologically mixing,
with slopes
$\pm s$,  with associated partition $a_0 < a_1 < \cdots < a_n$. Assume
\begin{enumerate}
\item $s$ is transcendental,
\item the lengths of the
intervals $[a_0,a_1], \ldots, [a_{n-2},a_{n-1}]$ are independent over $\Z[s,s^{-1}]$,
\item $\tau(\{0,1\}) \cap \{0,1\} = \emptyset$.
\end{enumerate}
\noindent Then the unique state on $DG(\tau)$ is an order isomorphism from  $DG(\tau)$ onto the subgroup $\sum_{i=1}^{n-1}
\Z[s,s^{-1}](a_i-a_{i-1})$ of $\R$, where the subgroup is given the order inherited from
$\R$. This isomorphism carries the action of $\L_*$ to multiplication by $s$.}

\prooff{Since $\tau$ is uniformly piecewise linear with slopes $\pm s$, then Lebesgue measure is scaled by $\tau$ by the factor $s$, 
and thus induces a state $\omega$ on $DG(\tau)$ which is scaled by $\L_*$ by the factor $s$ 
(Proposition \ref{1.56} and Proposition \ref{1.72}).
 As observed before Proposition \ref{1.96}, $\tau$ will be topologically exact. Therefore there
is a unique state
$\omega$ on $DG(\tau)$, which must coincide with the state given by Lebesgue measure, 
 and the order on $DG(\tau)$ is the strict order from that state, cf. Corollary \ref{1.75}. 

We are going to show that $\omega$ is a group and order isomorphism from $DG(\tau)$ onto  the subgroup
$\sum_{i=1}^{n-1} \Z[s,s^{-1}](a_i-a_{i-1})$ of $\R$.
We first show $\omega$ maps $DG(\tau)$ onto this subgroup. We begin by showing that $DG(\tau)$ is generated as a  $\Z[t,t^{-1}]$ module by $E_1,
E_2, \ldots, E_{n-1}$. Let $E_i = I(a_{i-1},a_i)$ for $1
\le i
\le n$. By Corollary \ref{1.42},
$E_1, E_2, \ldots, E_n$ generate the dimension module. (Recall that we identify intervals with the
equivalence classes in $DG(\tau)$ of their characteristic functions.) By the assumption (iii), if
$\tau(a_x) = 0$ and
$\tau(a_y) = 1$ with $x < y$, then $\sum_{i=x}^{y-1} \pm \L_*(I(a_i,a_{i+1})) = I(0,1)$, with appropriate
choices of the $\pm$ signs.  (See 
\cite[proof of Prop. 10.6]{Shu} for details.) Thus
$I(0,1)$ is in the submodule generated by
$E_1,
\ldots, E_{n-1}$.  Since $E_n = I(0,1) - \sum_{i=1}^{n-1} E_i$, then $E_n$ also is in the submodule generated
by $E_1, E_2, \ldots, E_{n-1}$,  so the latter elements generate
$DG(\tau)$.  

Thus every element of $DG(\tau)$ has the form $\sum_{i=1}^{n-1} p_i(\L_*)E_i$, where $p_1, \ldots, p_{n-1} \in \Z[t,t^{-1}]$. 
Since $\omega$ is scaled by $\L_*$ by the factor $s$, we have
$$\omega(\sum_{i=1}^{n-1} p_i(\L_*)E_i) = \sum_{i=1}^{n-1} p_i(s)\omega(E_i) = \sum_{i=1}^{n-1} p_i(s)(a_i-a_{i-1})\in
 \sum_{i=1}^{n-1}
 \Z[s,s^{-1}](a_i-a_{i-1}).$$ 
It follows that $\omega$ maps $DG(\tau)$ onto $\sum_{i=1}^{n-1}\Z[s,s^{-1}](a_i-a_{i-1})$. 

Next we show the kernel of $\omega$ is zero. 
Since
$s$ is transcendental, if $0 \not= p \in \Z[t,t^{-1}]$, then  $p(s)$ is non-zero. 
If $p_1, \ldots, p_{n-1} \in
\Z[t,t^{-1}]$ and
$g =
\sum_{i=1}^{n-1} p_i(\L_*)E_i$, then $\omega(g) =
\sum_{i=1}^{n-1} p_i(s) (a_i-a_{i-1})$, which by (ii) is non-zero
unless all $p_i(s)$ are zero. Since $s$ is transcendental, this can happen only if all $p_i$ are zero. Thus 
the kernel of $\omega$ is zero.

We've shown $\omega$ is a group isomorphism from $DG(\tau)$ onto $\sum_{i=1}^{n-1}\Z[s,s^{-1}](a_i-a_{i-1})$. 
Since $DG(\tau)$ has the strict order from $\omega$, then $\omega$ is also an order isomorphism. Finally, since $\omega$
is scaled by $\L_*$ by the factor $s$, by the definition of scaling (see above Proposition \ref{1.72}),
 $\omega$ carries the action of $\L_*$ to multiplication by $s$, which completes the proof.}

If $\tau:I\to I$ is continuous, mixing, and \pwm, then $\tau$ is conjugate to a uniformly piecewise linear map whose
slopes are $\pm s$ with $s = \exp h_\tau$, cf. Corollary \ref{1.65}. If $a_0, a_1, \ldots, a_n$ is the partition associated with
$\tau$, and
$\mu$ is the unique measure on $I$ scaled by $\tau$,
then the conjugacy can be chosen to carry $(a_{i-1},a_i)$ onto an interval of length $\mu([a_{i-1},a_i])$, cf. Proposition
\ref{1.59}. Thus Propositions
\ref{1.96} and \ref{1.97} are valid for such a map $\tau$ if $s$ is taken to be $\exp(h_\tau)$, and $(a_i-a_{i-1})$ is replaced by
$\mu([a_{i-1},a_i])$.

\exe{1.98}{Let $\tau:I\to I$ be the uniformly piecewise linear map in Figure \ref{fig1} on page~\pageref{fig1}, with slopes $\pm s$,
with $s$ transcendental, and $2 < s< 3$.  Since $s > 2$, under iteration, the length of any interval expands until the interval
contains both critical points, and then the next iterate equals $I$. Thus
$\tau$ is topologically exact.  Assume that the
first critical point is at
$\alpha
\notin
\Z[s,s^{-1}]$. Then by Proposition \ref{1.97}, the dimension module is order isomorphic to
$\alpha
\Z[s,s^{-1}] +
\Z[s,s^{-1}]$.  Now let $\tau'$ be the \pwm\ map whose graph is that of $\tau$ turned upside down, i.e., $\tau'(x) = 1-\tau(x)$.
Then the same argument shows that $\tau'$ is topologically exact, and that the dimension module for $\tau'$ also is order
isomorphic to
$\alpha
\Z[s,s^{-1}] +
\Z[s,s^{-1}]$.  Yet  $\tau$ and $\tau'$ are not conjugate. (For example, they have different numbers of fixed points.) Thus 
condition (i) in Theorem
\ref{1.79} is not redundant.}

\section{$\beta$-transformations}

If $\beta > 1$, the {\it $\beta$-transformation\/} on $[0,1)$ is the map $\tau_\beta: x
\mapsto
\beta x
\bmod 1$. (See Figure \ref{fig1} on page~\pageref{fig1}.) We extend $\tau_\beta$ to $[0,1]$ 
by defining $\tau_\beta(1)= \lim_{x\to
1^-}\,
\tau_\beta(x)$. Then $\tau_\beta$ is
\pwm, and the associated partition is $\{0, 1/\beta, 2/\beta, \ldots, n/\beta, 1\}$, where $n= [\beta]$ is the greatest integer
$\le \beta$.

Katayama, Matsumoto, and Watatani have associated C*-algebras $F_\beta^\infty$ and
$O_\beta$ with the $\beta$-shift. The algebra $F_\beta^\infty$ is a simple AF-algebra  (\cite[Prop. 3.5]{KatMatWat}), so
$K_0(F_\beta^\infty)$ is a simple dimension group.  In this section, we will compute the dimension group
for $\tau_\beta$, and will show that the dimension groups $DG(\tau_\beta)$ and
$K_0(F_\beta^\infty)$ are isomorphic as abelian groups, and are also order isomorphic if 1 is eventually periodic.

\lem{1.105}{For $\beta > 1$,   the $\beta$-transformation $\tau$ is topologically exact, and thus  $DG(\tau)$ is a simple
dimension group. There is a unique state, given by Lebesgue measure, and scaled by $\L_*$ by the factor $\beta$. The order on
$DG(\tau)$ is the strict order given by this state.}

\prooff{If $J$ is any interval whose interior does not contain one of the partition points $k/\beta$ for $1 \le k
\le n$,  (where $n$ is the greatest integer $\le \beta$), then applying $\tau$
multiplies the length of $J$ by $\beta$.  Thus for $n$ large enough, the interior of $\tau^n(J)$ includes
one of these partition points, and then the interior of $\htau^{n+1}(J)$ includes zero.  Succeeding images
will contain 0 as an interior point, and the component of the images containing 0 will  grow in length until it  includes
$[0,1]$. Thus
$\tau$ is exact.  Simplicity of $DG(\tau)$ follows from \cite[Cor. 5.4]{Shu}. Since Lebesgue measure is scaled by $\tau$ by the
factor $\beta$, then it induces a state,  scaled by $\L_*$ by the factor $\beta$, cf. Proposition \ref{1.72}. That there is a
unique state, and that  $DG(\tau)$ has the strict order from this state, follow from  exactness of $\tau$ and Corollary
\ref{1.75}.}

If $x \in I$, and $n=[\beta]$, we give the intervals
$[0,1/\beta)$, $[1/\beta,2/\beta)$,
$\ldots$, $[(n-1)/\beta, n/\beta)$, $[n/\beta, 1]$ the labels $0, 1, \ldots, n$.  (If $\beta = n$, there are $n+1$ intervals,
with the last one being the singleton $\{1\}$.)  We define the {\it itinerary\/} of $x\in I$ to be the sequence
$n_0n_1n_2\ldots$, where
$\tau^kx$ is in the interval with label $n_k$. (Warning: in general, this is not a Markov partition, so this is not an
itinerary map in the sense used in \cite{Shu}.)

If $\tau_\beta$ is Markov (cf. \S7), then by definition, the orbit of
each endpoint of each interval of monotonicity under $\hat
\tau_\beta$ must be finite,  and so in particular 1 has a finite
orbit, and hence is eventually periodic.  Conversely, if 1 is eventually
periodic, then the orbit of 1 together with the points 0,
$1/\beta$, $\ldots$, $n/\beta$ give a Markov partition. Thus the map $\tau_\beta$ will be Markov iff 1 is eventually periodic.
In that case, if $A$ is the associated incidence matrix, the dimension triple for $\tau_\beta$ is isomorphic to
$(G_A, G_A^+, A_*)$, cf. Definition \ref{1.81} and Proposition \ref{1.82}. In Proposition \ref{1.107}, we will express $G_A$ as
an inductive limit with respect to right multiplication on $\Z^q$ by a  matrix $B_\L$ whose entries  are expressed in terms of
the itinerary of 1, and will describe the unique state.

\lem{1.106}{Let $\tau$ be the $\beta$-transformation, and $\sigma:X\to X$ the associated
local homeomorphism. Assume
$1$ is eventually periodic, with $1, \tau 1, \ldots, \tau^{p-1} 1$
distinct, and with $\tau^p 1 = \tau^k 1$ for
some $k < p$.  Let 
the itinerary of 1 be
$n_0 n_1n_2
\ldots n_k \ldots n_{p}\ldots$, and let $M= \Z[t]I(0,1)$. Let $m(t)$ be the minimal
polynomial for $\L$ restricted to $M$. 
\begin{enumerate}
\item If $\tau 1 = 1$, then $m(t) = t-n_0$, and $M \cong \Z$. \label{1.106i}
\item If $\tau 1 \not= 1$ and \label{1.106ii}
$\tau^p 1
\not= 0$, then $M \cong \Z^p$, and
\begin{equation*}
m(t) = t^p - n_0 t^{p-1} - n_1t^{p-2}
- \cdots - n_{p-1}  - (t^k -
n_0t^{k-1} - \cdots - n_{k-1}).
\end{equation*}
\item If $\tau^p 1 = 0$, then $M \cong \Z^{p-1}$, and \label{1.106iii}
$$
m(t) = t^{p-1} - n_0 t^{p-2} - n_1t^{p-3}
- \cdots - n_{p-2}.$$
\end{enumerate}
}

\prooff{If $\tau 1 = 1$, then $\beta$ is an integer $n_0$, and $\L I(0,1) = n_0 I(0,1)$, so $m(t) = t-n_0$, and $M \cong \Z$ follow.
Hereafter we assume $\tau 1 \not= 1$.

Suppose that $\tau^p 1 \not= 0$, so that $\tau^{p-1}1 \not= 0$. The points
$0$, $1$, $\tau 1$, $ \tau^2 1$, $\ldots$, $\tau^{p-1} 1$ will be distinct.  Let $b_0 < b_1 < \cdots < b_p$ be these points
arranged in increasing order, and  define
$E_i = I(b_{i-1},b_i)$ for $1\le i \le p$. Then  every element of $M$ can be written
uniquely as an integral combination of $E_1, \ldots, E_p$, and $\sum_i z_i E_i \ge 0$
iff $(z_1, \ldots, z_p ) \ge 0$, so $M \cong \Z^p$. If instead $\tau^p 1 = 0$, then $\tau^{p-1}1 = 0$. Then there are $p$
distinct points among $0$, $1$, $\tau 1$, $ \tau^2 1$, $\ldots$, $\tau^{p-2} 1$, so a similar argument shows $M \cong \Z^{p-1}$.

Whether or not $\tau^p 1 = 0$, 
$\sigma$ maps each interval of monotonicity onto
$I(0,1)$ (except  the rightmost interval), so for $ 1\le j\le
p$,
\begin{equation}
\L I(0,\tau^{j-1}1) = n_{j-1}I(0,1) + I(0,\tau^{j} 1).
\end{equation}
Thus
\begin{equation}
I(0,\tau^{j}1)= \L I(0,\tau^{j-1}1) -
n_{j-1}I(0,1).\label{(1.36)}
\end{equation}
Hence, by induction,
\begin{equation}
I(0,\tau^j(1)) = (\L^j - n_0\L^{j-1} - n_1\L^{j-2} \cdots -
n_{j-1})I(0,1).\label{(1.37)}
\end{equation}

(ii)  Now we will establish the formula in (ii) for $m(t)$. Suppose that $\tau^p 1 \not= 0$. Since $\tau^k(1) = \tau^p(1)$, from
the versions of (\ref{(1.37)}) with
$j = k$ and $j = p$, it follows that $\L$ satisfies $m(\L)I(0,1) =
0$, where $m$ is the polynomial in (ii).
It is clear from (\ref{(1.37)}) that
the linear span in $C(X)$ of $I(0,1)$, $\L I(0,1)$, $\ldots$,
$\L^{p-1}I(0,1)$ is the same as that of
$I(0,1)$, $I(0,\tau 1)$, $\ldots$, $I(0,\tau^{p-1}1)$. 

Since the points $1$, $\tau
1$,
$\ldots$,
$\tau^{p-1}1$ are distinct, and none equal 0, then $I(0,1)$, $I(0,\tau 1)$,
$\ldots$, $I(0,\tau^{p-1} 1)$ are linearly independent in $C(X)$, so the
same must be true of the functions   $I(0,1)$, $\L I(0,1)$,
$\ldots$,
$\L^{p-1}I(0,1)$. It follows
that $\L|_M$ cannot satisfy a polynomial of degree less than $p$,
so $m$ as described in (\ref{1.106ii}) is the minimal polynomial for $\L$.

(iii) Assume $\tau^{p-1}1 = 0$.   By (\ref{(1.37)}),
$$(\L^{p-1} - n_0\L^{p-2} - n_1\L^{p-3} \cdots -
n_{p-2})I(0,1) = 0.$$
In this case, $I(0,1)$, $I(0,\tau 1)$,
$\ldots$, $I(0,\tau^{p-2} 1)$ are linearly independent in $C(X)$, so there is no polynomial of degree less than
$p-1$ that annihilates $\L|_M$. Thus the 
minimal polynomial for $\L|_M$ is as described in (\ref{1.106iii}).}

The notation used below for inductive limits was introduced after Definition \ref{1.81}.

\prop{1.107}{Let $\tau$ be the $\beta$-transformation, and assume that $1$ is eventually
periodic. Let  $M= \Z[t]I(0,1)$; cf. Lemma \ref{1.106}. Let the minimal polynomial $m$ of $\L$ on $M$ be
\begin{equation}
m(\lambda) =
\lambda^q - a_0 \lambda^{q-1} - a_1 \lambda^{q-2} - \cdots - a_{q-1}.
\end{equation}
Then
$DG(\tau)$  is isomorphic as a group to the stationary inductive limit given by right
multiplication on $\Z^q$ by the matrix
\begin{equation}
B_\L = 
\begin{pmatrix}0&1&0&\cdots&0\cr
        0&0&1&\cdots & 0\cr
       \cdots&\cdots&\cdots&\cdots&\cdots\cr
         0&0&0&\cdots & 1\cr
       a_0&a_1&a_2&\cdots& a_{q-1}\end{pmatrix}.
\end{equation}
The  unique state $\omega$ on $DG(\tau)$ is given by
Lebesgue measure.  On the inductive limit, $\omega$ satisfies 
\begin{equation}
\omega([((z_0, \ldots,z_{q-1}),n)]) = \beta^{-n}\sum_i z_i \beta^i.
\end{equation}
The order on $DG(\tau)$
is the strict order given by this state, and the range of this state is
$\Z[\beta,\beta^{-1}]$.}

\prooff{Note that $\L(M) \subset M$. We will prove that
\begin{enumerate}
\item for each $f \in C(X,\Z)$, there exists $n\ge 0$ such that  $\L^n f \in M$.
\item $\L_*$ is surjective.
\end{enumerate}
It will follow that $DG(\tau)$ is isomorphic to the stationary inductive limit $M\mapright{\L}M$, cf. \cite[Lemma 8.3]{Shu}.

Let $B = \{0, 1,\tau 1, \tau^2 1, \ldots,
\tau^{q-1}1\}$. Note that $\tau(B) \subset B$.  By (\ref{(1.37)}),  every element of $M$ is an integral combination of elements
$I(0,b)$ with
$b
\in B$, and thus is contained in the subgroup generated by $I(b_i,b_j)$ for $b_i, b_j \in B$. Let $C = \{0, 1/\beta, \ldots,
n/\beta, 1\}$, where
$n = [\beta]$. Note that
$\htau(C)
\subset B$. By definition of $I_1$, the orbit of every point in $I_1$ eventually lands in the orbit of
$C$, and therefore eventually lands in $B$. From this (i) follows. Since a Markov map is  eventually surjective, then $\L_*$ is
surjective, cf. \cite[Lemma 3.5]{Shu}, so (ii) holds.

 By Lemma \ref{1.106}, $I(0,1)$, $\L I(0,1)$, $\ldots$, $\L^{q-1} I(0,1)$ are a
basis for $M\cong \Z^q$, and with respect to this basis, $\L$ has the matrix $B_\L$
given above. Define $\psi:\Z^q\to C(X,\Z)$ by $\psi(z_0, z_1, \ldots, z_{q-1}) = \sum_i z_i \L^i I(0,1)$.  Then
$\psi$ is an  isomorphism from $\Z^q$ onto $M = \Z[t]I(0,1)$ (as abelian groups), and
satisfies
$$\psi(vB_\L) =
\L \psi(v)  \text{ for all $v \in \Z^q$.}$$
It follows that $\psi$  induces an isomorphism $\Phi$ of the inductive limit $(G_{B_\L},B_\L)$ and
the stationary inductive limit $M\mapright{\L}M$. Thus $(DG(\tau), \L_*) \cong  (G_{B_\L},B_\L)$.
(However, since $\psi$ is not positive, this is not an order isomorphism.)

 Recall that $DG(\tau)$ has a unique state, given by Lebesgue measure $\mu$, and has
 the strict order given by that state (Lemma \ref{1.105}).  
If we equip $(G_{B_\L},B_\L)$ with the order carried over by the isomorphism $\Phi$ from $DG(\tau)$,
then
$(G_{B_\L},B_\L)$ has the strict order given by the homomorphism $\omega:G_{B_\L} \to \R$ defined by $\omega 
=\mu\circ\Phi$. Since Lebesgue measure is scaled by $\tau$ by the factor $\beta$, then $\L_*$ scales the associated
state by the same factor.  Thus for
$v = (z_0, z_1,
\ldots, z_{q-1})\in \Z^q$,
$$\omega([v,n]) = \mu(\L_*^{-n}[\psi(v)]) = \beta^{-n} \sum_i z_i\mu(\L^i I(0,1))= \beta^{-n} \sum_i z_i\beta^i.$$}

We need the following result from \cite{Shu}. Recall that a module is \emph{cyclic} if it is singly generated.

\prop{1.47}{(\cite[Prop. 7.3]{Shu}) Let $\tau:I\to I$ be \pwm,  with associated partition $C$. Assume that $DG(\tau)$
is cyclic, and that there exists $a \in \{0,1\}$ with an infinite orbit, such that 
\begin{equation}
\htau (C\setminus \{a\}) \subset C.\label{(1.16)}
\end{equation}
 Then $DG(\tau)\cong \Z[t,t^{-1}]$ as abelian groups, with the action of $\L_*$ given by
multiplication by $t$.}

\prop{1.108}{Let $\tau$ be the $\beta$-transformation. If $1$ is not eventually
periodic, then
$DG(\tau)$ is  isomorphic as a group to $\Z[t,t^{-1}]$, with the action of $\L_*$ given by
multiplication by $t$. The unique state is given by $p\mapsto p(\beta)$, and $DG(\tau)$ has
the strict order from this state.}

\prooff{Each of the intervals
$I(i/\beta, (i+1)/\beta)$ is mapped onto
$I(0,1)$ by
$\L$, and the jumps at partition points are all of size 1.  From Theorem \ref{1.41}, it follows that $I(0,1)$ generates the
module
$DG(\tau)$. If $a$ is a partition point other than $1$, then $\htau(a) \subset \{0,1\}$, and 1 is not eventually periodic,
so
$DG(\tau)
\cong
\Z[t,t^{-1}]$ follows from Proposition \ref{1.47}.  Lebesgue measure scales $\tau$ by the factor $\beta$, so takes
$p(\L)I(0,1)$ to $p(\beta)$.  The remaining statements follow from Lemma \ref{1.105}).}

Note  in Proposition \ref{1.108} that when $1$ is not eventually periodic, the range of the unique state is
$\Z[\beta,\beta^{-1}]$. There will be infinitesimals iff $\beta$ is algebraic, cf. Proposition \ref{1.77}.  If $\beta$ is
transcendental, then $DG(\tau)$ is isomorphic as an ordered group to $\Z[\beta,\beta^{-1}]$, and the action of
$\L_*$ on the latter is given by multiplication by $\beta$.

 We now show that the dimension group $DG(F_\beta^\infty)$ defined in \cite{KatMatWat} is isomorphic to $DG(\tau_\beta)$.

\prop{1.109}{If $\beta> 1$ and $\tau$ is the $\beta$-transformation, then the dimension groups $DG(\tau)$ and
$K_0(F_\beta^\infty)$ are isomorphic as abelian groups, and are also order isomorphic if 1 is eventually periodic.}

\prooff{If $x \in [0,1]$, then   the
{\it $\beta$-expansion\/} of $x$ is $x = \sum_{i=1}^\infty \frac{\eta_i}{\beta^i}$, where for each $k$,
$\eta_k$ is the greatest integer in $\beta^k (x-\sum_{i=1}^{k-1} \frac{\eta_i}{\beta^i})$. The
itinerary of $1$ defined above Lemma \ref{1.106} will be the sequence of coefficients $\eta_1 \eta_2 \eta_3\cdots$.   If $1$ is
eventually periodic,  Lemma \ref{1.105} and \cite[Prop. 3.5 and Lemma 6.2]{KatMatWat} imply that $DG(\tau_\beta)$ and
$K_0(F_\beta^\infty)$ are simple dimension groups with unique states. Then   Proposition \ref{1.107} and \cite[Lemma
4.9]{KatMatWat} show these groups are given by the same inductive limit, and have the same  state, and thus the same
strict order. If 1 is not eventually periodic, the desired group isomorphism follows from Proposition
\ref{1.108} and \cite[Lemma 4.9, and proof of Thm. 6.1]{KatMatWat}.}

{\bf Remarks} It is an open question whether $DG(\tau)$ and
$K_0(F_\beta^\infty)$ are isomorphic as ordered abelian groups when $1$ is not eventually periodic.

Let $S:[0,1)\to \{0, 1, \ldots,n-1\}^\N$ be the itinerary map for the $\beta$-transformation described above. Let $X_\beta$ be
the closure of the image of $S$ in $\{0, 1, \ldots, n-1\}^\N$, where $n = [\beta]$, and let $\sigma_\beta$ be the shift on
$X_\beta$. Then
$\sigma_\beta$ is called the $\beta$-shift.  We note that $S(I_0)$ will be dense in $S([0,1))$, so $X_\beta =
\overline{S(I_0)}$.

Also denote by $S$ the corresponding itinerary map on $X$ with respect to the partition $I(0,1\beta)$, $I(1/\beta,2/\beta)$,
$\ldots$. Then the image of
$S$ on
$X_0$ is the same as on $I_0$, so $S(X_0)$ is dense in $X_\beta$. Since $X_0$ is dense in $X$ and $S$ is continuous on $X$, then
$S(X) = X_\beta$, and $S$ will be a semi-conjugacy from $(X,\sigma)$ onto $(X_\beta,\sigma_\beta)$.

If the orbit of 1 lands on 0, or $\tau 1 = 1$, then $\tau$ is Markov, and the itinerary map will be 1-1 on $X$, and thus will be
a conjugacy from $(X,\sigma)$ onto $(X_\beta, \sigma_\beta)$. Since
$\tau$ and
$\sigma$ are topologically exact, then $(X,\sigma)$ is a shift of finite type \cite[Prop. 8.5]{Shu}, and thus so is
$(X_\beta,\sigma_\beta)$.

 If the orbit of 1 is
finite but doesn't land on 0, and
$\tau 1 \not= 1$, again $\tau$ is Markov, but now the itinerary map $S$ with respect to the partition of $X$ associated with
$\{0, 1/\beta,
\ldots, n/\beta, 1\}$ will not be 1-1.  For example, if $x = \tau 1$, then $x^+$ and $x^-$ will have the same itineraries.  In
this case the itinerary map is just a semi-conjugacy from the shift of finite type $(X,\sigma)$ onto the sofic shift
$(X_\beta,\sigma_\beta)$.

\appendix 

\section{Convergence of powers of the Perron-Frobenius operator}

In this appendix  $\tau:I\to I$ will be a \pwm\ map, with associated local
homeomorphism $\sigma:X\to X$. We will investigate the convergence of the powers of the
Perron-Frobenius operator associated with $\sigma$. 
\medskip

{\bf Standing Assumption\/} Throughout this appendix, we will assume that there
exists a non-atomic probability measure $\mu$ on $X$, with full support, scaled by
$\sigma$ by a factor
$s > 1$ (cf. Definition \ref{1.54}).   
\medskip

 The existence of such a scaling measure is equivalent to $\tau$ being conjugate to a
piecewise linear map with slopes $\pm s$ with $s > 1$, cf. Propositions \ref{1.56} and
\ref{1.59}. By Proposition  \ref{1.64}, such a
measure will exist if
$\tau$ is transitive and not essentially injective.  Below, ``a.e." will mean with respect to the measure $\mu$.

We will begin with a summary of some results in a measure-theoretic context: convergence
results are in the space $\bv$ (defined below) of functions of bounded variation in $L^1$, equalities
of functions hold a.e., and partitions of $X$ are a.e. These
results are well known; mainly we will be following Rychlik \cite{Ry}. An exposition of convergence results for
the Perron-Frobenius operator can be found in \cite{BG}.

For our application,
it will be important to establish the corresponding results in a topological context, so that functions are continuous, 
convergence of sequences of functions is uniform,  equalities of functions hold everywhere,  and partitions of $X$
are into clopen sets. For this purpose the key is working with the Perron-Frobenius operator associated with
$\sigma$, rather than the one associated with
$\tau$. Since
$\sigma$ is a local homeomorphism, then the transfer operator $\L_\sigma$ (cf. Definition \ref{1.22})  maps $C(X)$ into $C(X)$,
which is what allows us to transfer results to the topological context. For related results for the special case of
$\beta$-transformations, see \cite{WalBeta}.

Let $P:L^1(X,\mu) \to
L^1(X,\mu) $ be the map given by $P= (1/s)\L_\sigma$.  Since $\sigma$ scales $\mu$ by the factor $s$,
then for $f \in \L^1(X,\mu)$, $\mu(\L f) = s \mu(f)$ (Lemma \ref{1.55}). Thus for $0 \le f$, 
$\|Pf\|_1 = \|f\|_1$, from which it follows that $\|P\|_1 =1$.  The map $P$ is the Perron-Frobenius operator
associated with
$\sigma$, i.e., for all $f \in L^1(X,\mu)$ and all Borel sets $A \subset X$, it satisfies
\begin{equation}
\int_A (Pf) \,d\mu =  \int_{\sigma^{-1}(A)} f\,d\mu
\end{equation}
(cf. \cite{Ry}).  

 For
$f:X\to \C$, recall that the {\it variation\/} of $f$ is
\begin{equation}
\var f = \sup\sum_{i=1}^n |f(a_i)-f(a_{i-1})|,
\end{equation}
where the supremum is over all finite sequences $a_0 < a_1 < \cdots
< a_n$ of points in $X$.  We say $f$ is of {\it bounded variation\/}
if $\var f < \infty$. Every function $f$ of bounded variation on $X$ is the
difference of increasing functions. (In fact, if
$a$ is the leftmost point of $X$, then we can write $f = g - h$ where
$g(x) = \var_{[a,x] \cap X}f$ and $h = g - f$, cf. \cite[\S1.4]{Riesz-Nagy}). It
follows that
$f$ has left and right limits at each point, and that $f$ is continuous except
at a countable set of points.

\defi{1.113}{$\bv$ is the space of functions in $L^1(X,\mu)$ equal a.e. to a
function of bounded variation, with the norm
\begin{equation}
\|f\|_\bv = \inf\{\,\var g \mid g = f \ {\rm a.e.}\} + \|f\|_1
\end{equation}
where $\|f\|_1$ denotes the norm from $L^1(X,\mu)$.}

Note that $\bv$ is a Banach space, cf. \cite{Ry}. Each function of
bounded variation is bounded, so
$\bv
\subset L^\infty(X,\mu)$. If $f \in \bv$ we write $\|f\|_\infty$ for the norm of
$f$ as a function in $L^\infty(X,\mu)$.

\lem{1.114}{If $f \in \bv$, then
\begin{equation}
\|f\|_\infty  \le \|f\|_\bv.
\end{equation}}

\prooff{Let $f$ have bounded variation. Since $\|f\|_1 \ge \inf_{x \in X}|f(x)|$, for each $\epsilon > 0$, there
is some
$y\in X$ such that $|f(y)| \le \|f\|_1+ \epsilon$. Then for each $x \in X$,
\begin{equation}
|f(x)| \le |f(x) - f(y)| + |f(y)| \le \var f + \|f\|_1+\epsilon,
\end{equation}
so $\|f\|_\infty \le\var f + \|f\|_1$.  As $f$ varies in its equivalence class in $L^1$, $\|f\|_\infty$ doesn't
change, and the infimum of the numbers $\var f + \|f\|_1$ is $\|f\|_\bv$. The inequality in the lemma
follows.}

\prop{1.115}{(Rychlik \cite[Thm. 1]{Ry}) $P$ maps $\bv$ into $\bv$. If $P$ is considered as an operator on
$\bv$, the set
$\lambda_1$,
$\lambda_2$,
$\ldots$,
$\lambda_q$ of points in the spectrum of $P$ with modulus 1 is finite and non-empty.    These points are 
all roots of unity,  and are  simple poles of the resolvent. There is a number
$r< 1$ such that all other points in the spectrum are contained in the disk around $0$ of
radius $r$. If for each $i$,  $Q_i$ is the
spectral projection corresponding to $\lambda_i$, then 
$P = \sum_i \lambda_i Q_i + R$,
where each $Q_i$ is a finite rank projection whose range is the eigenspace of $P$ for eigenvalue $\lambda_i$,
with $Q_iQ_j = 0$ for $i \not = j$, $Q_i R = RQ_i = 0$ for each $i$,  and
with the spectral radius of
$R$  less than 1.}

Let $N$ be the least positive
integer in Proposition \ref{1.115} such that $\lambda_i^N = 1$ for all $i$, and let $Q = \sum_i Q_i$. Note that $Q$
is a projection onto the space spanned by all eigenvectors corresponding to eigenvalues of modulus 1, or
equivalently, onto the space of fixed points of $P^N$.

\lem{1.116}{If $f \in C(X)\cap \bv$, then  $f$ has bounded variation.}

\prooff{Choose $g$ of bounded variation, such that $f = g$ a.e. Let $x_1< x_2 < \ldots < x_k$ be any points in
$X$. Since
$\mu$ has full support, in any open interval around a point $x \in X$, there exists a point $t$ such that $g(t)
= f(t)$. Given
$\epsilon > 0$,   choose points $t_1< t_2 < \ldots < t_k$ such that $\sum_i|f(x_i)-f(t_i)| < \epsilon/2k$
and such that $g(t_i) = f(t_i)$ for all $i$. Then 
\begin{equation}
\sum_i |g(t_i) - g(t_{i+1})| = \sum_i |f(t_i) - f(t_{i+1})| \ge \sum_i |f(x_i) - f(x_{i+1})|- \epsilon.
\end{equation}
It follows that $\var g \ge \var f$, so $f$ has bounded variation.}

\cor{1.117}{For all $f\in
\bv$, 
\begin{equation}
Qf = \lim_k P^{Nk}f,
\end{equation}
where convergence is in $\bv$.  If
$f \in C(X) \cap \bv$, then $P^{Nk}f$ converges uniformly to a function in $C(X)\cap \bv$,
and so $Qf$ is equal a.e. to a function in $C(X)\cap \bv$.}

\prooff{The first statement is an immediate consequence of Proposition \ref{1.115}. 
Recall that $P = (1/s)\L_\sigma$ maps $C(X)$ into $C(X)$ (as noted after Lemma \ref{1.55}). If $f \in C(X) \cap \bv$,
since the $\bv$ norm dominates the $L^\infty$ norm (Lemma \ref{1.114}),  then
$\|P^{Nk}f - Qf\|_\infty \to 0$. Hence $\{P^{Nk}f\}_{k=1}^\infty $ is a sequence of
continuous functions, and is Cauchy with respect to the supremum norm. Let $g$ be
the uniform limit of this sequence. Then $g$ is continuous, and $g = Qf$ a.e. By Lemma \ref{1.116}, $g$ has bounded
variation.}

 If $A$, $B$ are subsets of $X$, we write $A = B$ a.e. if the symmetric
difference
$(A\setminus B) \cup (B \setminus A)$ has measure zero with respect to $\mu$.  Note
that all equalities of functions in Theorem \ref{1.118} are a.e. In Theorem \ref{1.118}, $Q$ denotes the map  defined in
Proposition \ref{1.115}, and  $N$ is the integer defined in the remarks following that proposition.

\theo{1.118}{(Rychlik  \cite[Thms. 3, 4]{Ry}) 
Then there are nonnegative
functions $\phi_1, \dots, \phi_q \in \bv$ and $\psi_1, \ldots, \psi_q \in
L^\infty$ such that
\begin{enumerate}
\item  $Q$ has a unique extension to a bounded operator on $L^1$, and $Q$
maps
$L^1$ into $\bv$. \label{1.118i}
\item For every $f \in L^1$,  \label{1.118ii}
$$
Qf = \sum_{i=1}^q \mu(\psi_i  f) \phi_i.
$$
\item  \label{1.118iii} $P^N\phi_i = \phi_i$, $\psi_i\circ \sigma^N = \psi_i$ for $i =
1,
\ldots, q$.
\item  \label{1.118iv} $\mu(\phi_i\psi_j) = \delta_{ij}$, $\min(\psi_i, \psi_j) = 0=
\min(\phi_i, \phi_j)$ for $i\not= j$, $\mu(\phi_i) = 1$ for $i = 1, \ldots, q$. 
\item  \label{1.118v} There exist measurable sets $C_1, \ldots, C_q \subset X$ such
that $\psi_i = \chi_{C_i}$ a.e. for $i = 1, \ldots, q$, and $X =
\cup_{i=1}^q C_i$ a.e.
\item  \label{1.118vi} There exists a permutation $\omega$ of the set $\{1, \ldots, q\}$
such that
\begin{equation}
P\phi_i = \phi_{\omega(i)}, \qquad \psi_{\omega(i)}\circ \sigma = \psi_i
\quad \hbox{for $i = 1, \ldots, q$}, \label{eq39}
\end{equation}
and such that $\omega^N$ is the identity.  If $\omega$ is a cycle, then $q = N$.\label{eq41}
\end{enumerate}}

\prooff{This result is contained in \cite[Thms. 3, 4]{Ry}, except for the last statement in (\ref{1.118vi}).  Suppose $\omega$ is
a
$q$-cycle. 
 By (\ref{1.118iv}), $\phi_1$, $\phi_2$,
$\ldots$, $\phi_q$ are distinct, so by (\ref{eq41}),
$\phi_1$,
$P\phi_1$,
$\ldots$,
$P^{q-1}\phi_1$ are distinct. By (\ref{eq39}), $P^N\phi_1 = \phi_1$, which implies that $N \ge q$. 

Since $\omega$ is a $q$-cycle, $P^q \phi_i = \phi_i$ for all $i$. By (\ref{1.118ii}), the range of $Q$ is
spanned by $\phi_1$, $\ldots$, $\phi_q$, so $P^q Q = Q$. Let $\lambda$ be an eigenvalue of modulus
1 for $P$, with eigenvector $\psi$. By construction, $Q$ fixes each such eigenvector, so $P^q Q\psi = Q\psi$ implies that $P^q
\psi = \lambda^q \psi = \psi$. Hence $\lambda^q = 1$ for each eigenvalue of $P$ of modulus 1, so by the definition of $N$, $q
\ge N$. Thus $N = q$.}

By definition, the support of a function $f\in C(X)$ is the closure of the set $\{x \mid f(x) \not= 0\}$.
If
$f
\in C(X,\Z)$, then $\{x \mid f(x) \not= 0\}$ is already closed, so $\supp f$ consists of the points where $f$ is not
zero. From the definition of
$\L_\sigma$,
\begin{equation}
0 \le f \in C(X) \implies \supp(\L_\sigma f) = \sigma(\supp f).\label{(1.46)}
\end{equation}

We are now going to strengthen the conclusions of Theorem \ref{1.118} by showing that each
$\phi_i$ can be chosen to be  continuous, with clopen support. Recall that $\mu$ is a
measure of full support on $X$, scaled by $\sigma$ by the factor $s$.

\cor{1.119}{Each $\,\phi_i\,$ in Theorem \ref{1.118} can be chosen to be continuous and of bounded variation. For such a choice of $\phi_i$, 
if $f \in C(X) \cap \bv$ with $\supp f \subset \{x \mid \phi_i(x) > 0\}$,
then
$P^{Nk}f$ converges uniformly to
$\mu(f)\phi_i$.}

\prooff{By definition of $\bv$, each $\phi_i$ is equal a.e. to a
function of bounded variation, so without loss of generality we may assume $\phi_i$ has bounded
variation on
$X$. Then $\phi_i$ is continuous except on a countable set of points (and thus
a.e., since $\mu$ is  non-atomic). If
$\phi_i$ were zero at each point of continuity, then it would be zero a.e., which
would contradict $\mu(\psi_i\phi_i) = 1$. Thus there is an open subset $V$ of $X$ 
such that $\phi_i > 0$ on $V$. 

Now let $f$ be any non-negative continuous function of bounded variation, not identically zero, such that the support of $f$  is
contained in $\{x \mid \phi_i(x) > 0\}$. By Theorem \ref{1.118} (\ref{1.118iv}) and (\ref{1.118v}), $\{x \mid \phi_i(x) > 0\}
\subset C_i$ a.e., so the support of $f$ is contained in $C_i$ a.e. By  Theorem \ref{1.118} (\ref{1.118ii}),  $Qf =
\mu(f)\phi_i$, and since
$\mu$ has full support,
$\mu(f) > 0$.  

By Proposition \ref{1.117}, $P^{Nk}f$
converges uniformly to $\mu(f)\phi_i'$, where $\phi_i'$ in $C(X)\cap \bv$ is equal a.e. to $Qf$.  
Thus $\phi_i= \phi_i'$ a.e., so by replacing $\phi_i$ by $\phi_i'$, we may choose $\phi_i$ to be continuous.

We have proven that $P^{Nk}f$ converges uniformly to $\mu(f)\phi_i$ for $0 \le f \in C(X)\cap \bv$ with $\supp f \subset \{x
\mid \phi_i(x) > 0\}$.  If $f$ is not necessarily non-negative, the same convergence result follows by writing $f = f^+ - f^-$,
with
$f^+ = \max(f,0)$ and $f^- = -\min(f,0)$.}

Hereafter we assume each $\phi_i$ in Theorem \ref{1.118} has been chosen to be continuous and of bounded variation. Note that
with such a choice for each $\phi_i$,  the equalities $P^N \phi_i = \phi_i$ and $P\phi_i = \phi_{\omega(i)}$ in Theorem
\ref{1.118} (vi) are equalities of continuous functions, so will hold everywhere, not just a.e.

\prop{1.120}{For each $i$, $\{x \mid \phi_i(x) > 0\}$ is clopen (and thus equals the support of $\phi_i$).}

\prooff{Let $g:I_0\to \R$ be defined by $g(x) = \phi_i(\pi^{-1}(x))$, where $\pi:X\to I$ is the collapse map.  Since $\phi_i\in
C(X)\cap
\bv$, then
$\phi_i$ is of bounded variation (Lemma \ref{1.116}), so $g$ will be the difference of increasing
functions on $I_0$. It follows that $g$ can be extended to the difference of increasing
functions on $I$, i.e., to a function of bounded variation on $I$, which we also
label as $g$.

Let $m$ be the measure on $I$ scaled by $\tau$ that corresponds to $\mu$, cf.
Proposition \ref{1.56}.
 Recall
that there is a conjugacy of
$(I,\tau)$ onto a piecewise linear map with slopes
$\pm s$, which carries $m$ to Lebesgue measure (Proposition \ref{1.59}).  Therefore we may
assume, without loss of generality, that $\tau$ is piecewise linear with slopes $\pm
s$, and that $m$ is Lebesgue measure.

Now let $V= \{x \mid \phi_i(x) > 0\}$. By continuity of $\phi_i$, $V$ is open.  Let $W = \pi(V\cap
X_0)$.  Since $\pi:X_0\to I_0$ is a homeomorphism (Proposition \ref{1.2}), then $W$ is an open subset of $I_0$ on which $g$
is positive.  

Since $\tau'(x) = 1/s$ a.e, then $(1/s)L_\tau$ is the Perron-Frobenius operator $P_\tau: L^1(I,m) \to L^1(I,m)$,
cf. \cite[eqn. (4.3.4)]{BG}. Since $((1/s)\L_\sigma)^N  \phi_i =P_\sigma^N \phi_i  = \phi_i$, then
$P_\tau^N g = ((1/s)\L_\tau)^N g = g$ on $I_0$, so $(P_\tau)^N g = g $ a.e. Since $(P_\tau)^N = P_{\tau^N}$,
then  $P_{\tau^N}g = g$ a.e.

Choose $g_1:I\to I$ lower semi-continuous so that $g = g_1$ a.e., 
cf. \cite[Lemma 8.1.1]{BG}.   
Then $g_1$ is a fixed point of
the Perron-Frobenius operator $P_{\tau^N}$. Since $g_1$ is lower semi-continuous, by \cite[Thm. 8.2.3]{BG} there is a
constant
$\alpha > 0$  such that $g_1(x)\ge\alpha$ for all $x$ such that $g_1(x) > 0$. 
In particular, $g_1(x)\ge \alpha$ a.e. on $W$.  It follows that $\phi_i(x) \ge \alpha$ a.e. on
the set $V\cap X_0$, which is dense in $V$.    By  continuity of $\phi_i$, and density of $X_0$,
$\phi_i\ge\alpha$ on
$V$. Thus
$V =
\{x
\in X
\mid \phi_i(x) \ge
\alpha\}$, so
$V$ is clopen.}

\cor{1.121}{With the notation of Theorem \ref{1.118}, let
$X_i= \supp \phi_i$ for $1 \le i \le q$. For each $i$, $\sigma(X_i) =X_{\omega(i)}$, and
$\sigma^N$ leaves $X_i$ invariant and is topologically exact on $X_i$. The sets $X_1$, $\ldots$, $X_q$ are disjoint.}

\prooff{Since $\min(\phi_i,\phi_j)= 0$ for $i \not= j$,
then $X_i\cap X_j =
\emptyset$ for $i \not= j$. By (\ref{(1.46)}), the support of $\phi_{\omega(i)} = (1/s)\L_\sigma\phi_i$ is
\begin{equation}
\supp (\L_\sigma \phi_i) = \sigma(\supp \phi_i) =
\sigma(X_i),
\end{equation}
 so
$\sigma(X_i) = X_{\omega(i)}$.   Since $\omega^N$ is the identity permutation, then $\sigma^N X_i = X_i$ for all $i$.

 Fix $i$ and let $V$
be any open subset of $X_i$.
Now choose $\epsilon > 0$ so that $\phi_i > \epsilon$ on $X_i$. Choose a 
nonnegative continuous function $f$ of bounded variation, not identically zero,  such that $f$ is zero
outside
$V$. By Corollary \ref{1.119}, $s^{-Nk} \L^{Nk}f$ converges
to $\alpha \phi_i$, where $\alpha = \mu(f) \not= 0$.

Choose
$k > 0$ such that
$\|s^{-Nk} \L^{Nk}f - \alpha \phi_i\|_\infty < \alpha\epsilon/2$. It follows that
$\L^{Nk}f> 0$ on all of $X_i$. Since $\L^{Nk}f$ is zero outside
$\sigma^{Nk}V$, then $\sigma^{Nk}(V) = X_i$. This completes the proof that $\sigma^N$ is
exact on $X_i$.}

{\bf Remarks} (i) Recall that our standing assumption in this appendix requires the existence of a measure scaled by $\tau$ by a
factor $s > 1$. If $s =1 $, then the results above can fail. (See the example after Corollary \ref{1.68}.)

(ii) Simple examples show that the union of the sets $X_i$ need not be all of $X$. However, this will be the case if
$\tau$ is transitive, cf. Theorem \ref{1.66}.

\end{document}